\documentclass[11pt]{article}

\usepackage{epsfig,amssymb,latexsym}
\usepackage{amsfonts,psfrag,amsmath,bbm,color}
\usepackage{cancel}

\usepackage{mathrsfs}
\usepackage{graphicx}
\usepackage{textcomp}
\usepackage{multirow}
\usepackage{enumerate}
\usepackage{cancel}
\usepackage{caption}  

\newcommand{\cE}{\mathcal{E}}

\newcommand{\Om}{\Omega}

\newtheorem{remark}{Remark}[section]
\newtheorem{assumption}{Assumption}[section]
\newtheorem{lemma}{Lemma}[section]

\newtheorem{definition}{Definition}[section]

\newcommand{\by}{{\bf y}}

\newcommand{\bff}{{\bf f}}

\def\PP{{{\rm l}\kern - .15em {\rm P} }}
\def\PN2{{\PP_{N}-\PP_{N-2}}}

\newcommand{\R}{\mathbb{R}}

\newcommand\figcaption{\def\@captype{figure}\caption}
\newcommand\tabcaption{\def\@captype{table}\caption}

\graphicspath{{./Figures/}}

\title{Are the Snapshot Difference Quotients Needed in the Proper Orthogonal Decomposition?}

\author{Traian Iliescu 
		\thanks{Department of Mathematics,
				Virginia Tech,
				456 McBryde Hall,
				Blacksburg, Va 24061
				({\tt iliescu@vt.edu}).}
        \and Zhu Wang
        	\thanks{Institute for Mathematics and its Applications, University of Minnesota, 354 Lind Hall, Minneapolis, MN 55455 ({\tt wangzhu@ima.umn.edu}). 
	}
	}

\begin{document}

\maketitle

\begin{abstract}
This paper presents a theoretical and numerical investigation of the following practical question:
Should the time difference quotients of the snapshots be used to generate the proper orthogonal decomposition basis functions?
The answer to this question is important, since some published numerical studies use the time difference quotients, whereas other numerical studies do not.
The criterion used in this paper to answer this question is the rate of convergence of the error of the reduced order model with respect to the number of proper orthogonal decomposition basis functions.
Two cases are considered: the $no\_DQ$ case, in which the snapshot difference quotients are not used, and the $DQ$ case, in which the snapshot difference quotients are used.
The error estimates suggest that the convergence rates in the $C^0(L^2)$-norm and in the $C^0(H^1)$-norm are optimal for the $DQ$ case, but suboptimal for the $no\_DQ$ case.
The convergence rates in the $L^2(H^1)$-norm are optimal for both the $DQ$ case and the $no\_DQ$ case.
Numerical tests are conducted on the heat equation and on the Burgers equation.
The numerical results support the conclusions drawn from the theoretical error estimates.
Overall, the theoretical and numerical results strongly suggest that, in order to achieve optimal pointwise in time rates of convergence with respect to the number of proper orthogonal decomposition basis functions, one should use the snapshot difference quotients.

\end{abstract}

\smallskip
\noindent \textbf{Keywords:}
proper orthogonal decomposition,
reduced order modeling,
error analysis


\pagestyle{myheadings}
\thispagestyle{plain}
\markboth{TRAIAN ILIESCU AND ZHU WANG}{{\scriptsize DIFFERENCE QUOTIENTS IN PROPER ORTHOGONAL DECOMPOSITION}}

\section{Introduction}
\label{s_introduction}

This paper addresses the following question:
{\it ``Should the time difference quotients (DQs) of the snapshots be used in the generation of the Proper Orthogonal Decomposition (POD) basis functions?"}

We emphasize that this is an important question.
There are two schools of thought: one uses the DQs (see, e.g. \cite{KV01,KV02,chaturantabut2010nonlinear,chaturantabut2012state}), the other does not (see, e.g. \cite{LCNY08,luo2009finite,chapelle2012galerkin,singler2012new}).

To our knowledge, the first instance in which the snapshot DQs were incorporated in the generation of the POD basis functions was the pioneering paper of Kunisch and Volkwein \cite{KV01}.
In that report, the authors considered two types of errors for a general parabolic equation: the time discretization errors and the POD truncation errors.  
They argued that one needs to include the temporal difference quotients in the set of snapshots; otherwise, the error will be suboptimal, containing an extra $\frac{1}{\Delta t^2}$ factor (see Remark 1 in \cite{KV01}).
Thus, the motivation for using the temporal difference quotients was purely theoretical.
In numerical investigations, however, the authors reported contradictory findings: in \cite{KV01}, the use of the DQs did not improve the quality of the reduced-order model; in \cite{homberg2003control}, however, it did.
Kunisch and Volkwein used again the snapshot DQs when they considered the Navier-Stokes equations \cite{KV02}.

The snapshot DQs were also used in the {\it Discrete Empirical Interpolation Method (DEIM)} of Chaturantabut and Sorensen \cite{chaturantabut2010nonlinear,chaturantabut2012state} (which is a discrete variant of the {\it Empirical Interpolation Method (EIM)} \cite{barrault2004eim}).
The motivation in \cite{chaturantabut2010nonlinear,chaturantabut2012state}, however, was different from that in \cite{KV01}.
Indeed, the authors considered in \cite{chaturantabut2010nonlinear,chaturantabut2012state} a general, nonlinear system of equations of the form $\by' = \bff(\by,t)$.
In the set of snapshots, they included not only the state variables $\by$, but also the {\it nonlinear} snapshots $\bff(\by,t)$.
They further noted (see page 48 in \cite{chaturantabut2012state}) that, since $\bff(\by,t) = \by'$ and $( \by^{n+1} - \by^{n} ) / \Delta t \sim \by'$, this is similar to including the temporal difference quotients, as done in \cite{KV01,KV02}.

To our knowledge, the first reports on POD analysis in which the DQs were not used were \cite{LCNY08} for the heat equation and \cite{luo2009finite} for the Navier-Stokes equations.
Chapelle et al.~\cite{chapelle2012galerkin} used a different approach that did not use the DQs either.
This approach employed the $L^2$ projection instead of the standard $H^1$ projection used in, e.g., \cite{KV01,KV02}.
Further improvements to the approach used in~\cite{chapelle2012galerkin}  (as well as that used in \cite{KV01,KV02}) were made by Singler in \cite{singler2012new}.

From the above discussion, it is clear that the question whether the snapshot DQs should be included or not in the set of snapshots is important.
To our knowledge, this question is still open.
This report represents a first step in answering this question.
All our discussion will be centered around the heat equation, although most (if not all) of it could be extended to general parabolic equations in a straightforward manner.

From a theoretical point of view, the only motivation for using the snapshot DQs was given in Remark 1 in \cite{KV01}.
The main point of this remark is the following:
In the error analysis of the evolution equation, to approximate $u_t(t^n)$, the time derivative of the exact solution $u$ evaluated at time $t^n$, one usually uses the DQ $\overline{\partial} u(t^n) := \frac{u(t^{n+1}) - u(t^{n})}{\Delta t}$.
To approximate the DQ $\overline{\partial} u(t^n)$ in the POD space, one naturally uses the POD DQ $\overline{\partial} u_r(t^n) := \frac{u_r(t^{n+1}) - u_r(t^{n})}{\Delta t}$, where $u_r$ is the POD reduced order model approximation.
We assume that $u_r(t^{n+1})$ is an optimal approximation for $u(t^{n+1})$ and that $u_r(t^{n})$ is an optimal approximation for $u(t^{n})$, where the optimality is with respect to $r$ and $\Delta t$.
Then, it would appear that, with respect to $\Delta t$, $\overline{\partial} u_r(t^n)$ is a {\it suboptimal} approximation for $\overline{\partial} u(t^n)$, because of the $\Delta t$ in the denominator of the two difference quotients.

Although the argument above, used in Remark 1 in \cite{KV01} to motivate the inclusion of the snapshot DQs in the derivation of the POD basis, seems natural, we point out that this issue should be treated more carefully.
Indeed, in the finite element approximation of parabolic equations, it is well known that the DQs $\overline{\partial} u_h(t^n) := \frac{u_h(t^{n+1}) - u_h(t^{n})}{\Delta t}$ are actually {\it optimal} (with respect to $\Delta t$) approximations of the DQs $\overline{\partial} u(t^n)$ (see, e.g., \cite{labovschii2009defect,shan2012decoupling}).
Thus, since the POD and finite element approximations are similar (both use a Galerkin projection in the spatial discretization), one could question the validity of the argument used in Remark 1 in \cite{KV01}.
We emphasize that we are not claiming that the above argument is not valid in a POD setting; we are merely pointing out that a rigorous numerical analysis is needed before drawing any conclusions.

\begin{figure}[htp]
\begin{center}
\begin{minipage}[h]{.6\linewidth}
\includegraphics[width=0.9\textwidth]{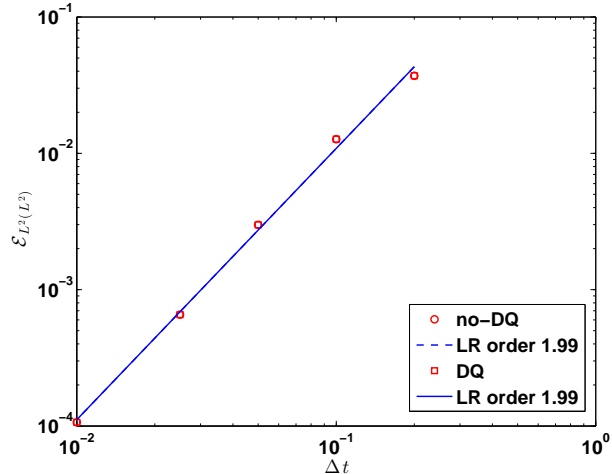}
\end{minipage}
\end{center}
\caption{
	Heat equation.
	Plots of the errors in the $L^2(L^2)$-norm with respect to the time step $\Delta t$ when the DQs are used (denoted by DQ) and when the DQs are not used (denoted by no-DQ).
	}
\label{fig:test1_dt}
\end{figure}

Our preliminary numerical studies indicate that not using the DQs does {\it not} yield suboptimal (with respect to $\Delta t$) error estimates.
For the heat equation (see Section~\ref{sec:numerical} for details regarding the numerical simulation), we monitor the rates of convergence with respect to $\Delta t$ for the POD reduced order model.
We consider two cases: when the DQs are used in the generation of the POD basis (the corresponding results are denoted by DQ), and when the DQs are not used in the generation of the POD basis (the corresponding results are denoted by no-DQ).
The errors (defined in Section~\ref{sec:numerical}) are listed in Table \ref{tab:test1_dt} and plotted in Figure \ref{fig:test1_dt} with associated {\it linear regressions (LR)}. 
Both no-DQ and DQ approaches yield an optimal approximation order $\mathcal{O}(\Delta t^2)$ in the $L^2$-norm. 

\begin{table}[h]
\centering
\tabcaption{Errors of the $no\_DQ$ and $DQ$ approaches when $\Delta t$ varies.}
\label{tab:test1_dt}
\begin{tabular}{c|ccc|ccc}
\hline
\multirow{2}{*}{$\Delta t$}&\multicolumn{3}{c|}{$no\_DQ$}&\multicolumn{3}{c}{$DQ$}\\
\cline{2-7}
{} & r  & $\cE_{L^2(L^2)}$  & $\cE_{L^2(H_1)}$ &  r  & $\cE_{L^2(L^2)}$  & $\cE_{L^2(H_1)}$  \\
 \hline
 2.00e-01 & 6  & 3.71e-02    &   9.26e-01   &  6  & 3.71e-02   &   9.26e-01      \\
 1.00e-01 & 11  & 1.27e-02   &   5.81e-01   &  11  & 1.27e-02   &   5.81e-01      \\
 5.00e-02 & 21  & 2.99e-03   &   1.97e-01   & 21  & 2.99e-03   &   1.97e-01      \\
 2.50e-02 & 41  & 6.53e-04   &   3.81e-02   & 41  & 6.53e-04   &   3.81e-02      \\
 1.00e-02 & 59  & 1.03e-04   &   1.15e-02   & 88  & 1.03e-04   &   1.15e-02   \\
 \hline
 \end{tabular}

 \end{table}

%
%
%

The rest of the paper is organized as follows:
In Section~\ref{sec:pod}, we sketch the derivation of the POD reduced order model.
In Section~\ref{sec:error}, we carefully derive the error estimates for the POD reduced order model. 
We focus on the rates of convergence with respect to $r$, the number of POD basis functions.
In Section~\ref{sec:numerical}, we present numerical results for two test problems: the heat equation and the Burgers equation.
Finally, in Section~\ref{sec:conclusions}, we draw several conclusions.

\section{Proper Orthogonal Decomposition Reduced Order Modeling}
\label{sec:pod}

In this section, we sketch the derivation of the standard POD Galerkin reduced order model for the heat equation.
For a detailed presentation of reduced order modeling in general settings, the reader is referred to, e.g.,  \cite{HLB96,KV99,Sir87abc,burkardt2006pod,AK04,bui2008model,wang2012proper,balajewicz2012novel}.

\begin{figure}[h]
\begin{center}
\begin{minipage}[h]{.9\linewidth}
\includegraphics[width=0.9\textwidth]{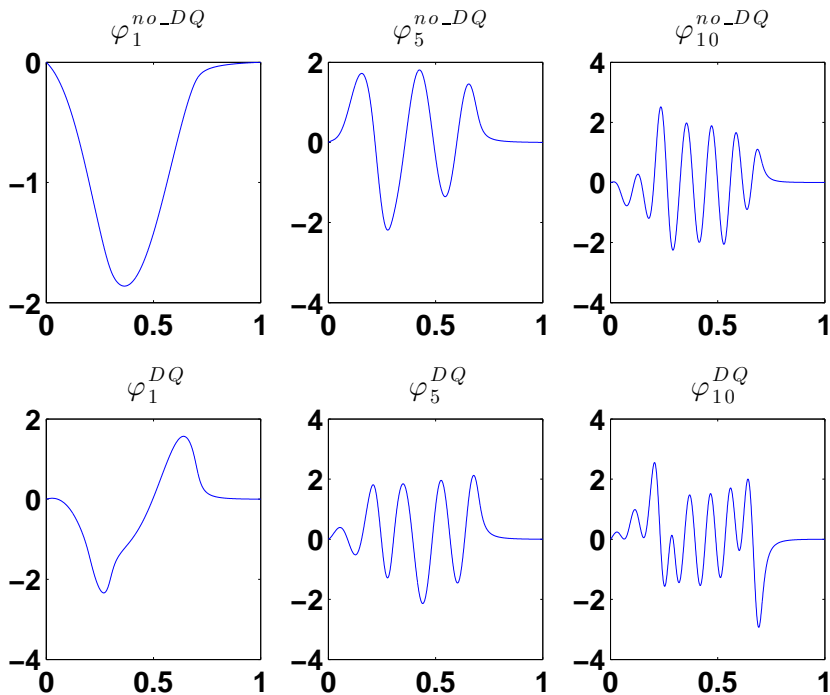}
\end{minipage}
\end{center}
\caption{
	Heat equation.
	Plots of the POD basis functions when the DQs are used (denoted by DQ) and when the DQs are not used (denoted by no-DQ).
	}
\label{fig:pod_basis_heat}
\end{figure}

\begin{figure}[h]
\begin{center}
\begin{minipage}[h]{.9\linewidth}
\includegraphics[width=0.9\textwidth]{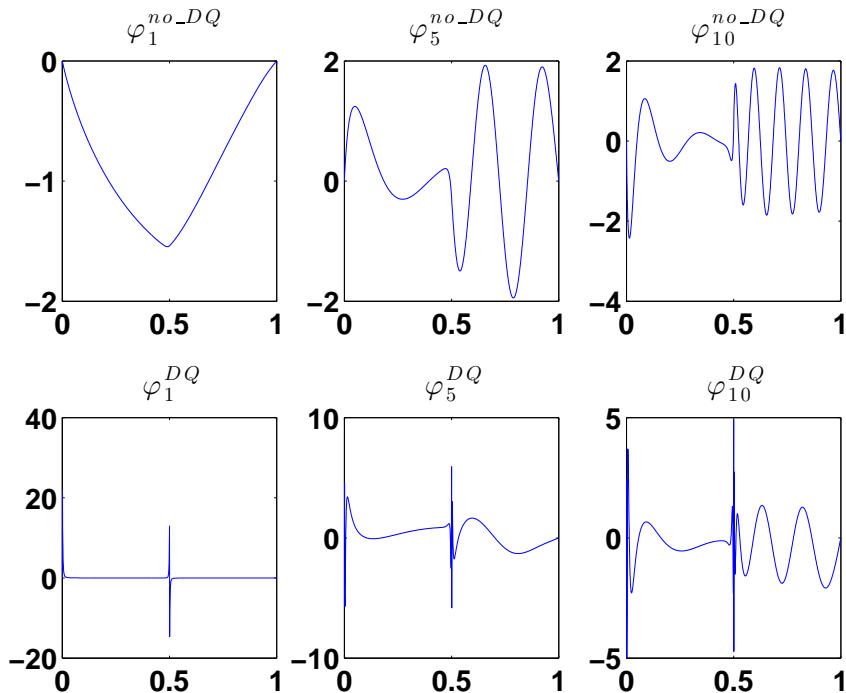}
\end{minipage}
\end{center}
\caption{
	Burgers equation.
	Plots of the POD basis functions when the DQs are used (denoted by DQ) and when the DQs are not used (denoted by no-DQ).
	}
\label{fig:pod_basis_burgers}
\end{figure}

For clarity, we will denote by $C$ a generic constant that can depend on all the parameters in the system, except on 
$M$ (the number of snapshots),
$d$ (the dimension of the set of snapshots, $V$), and
$r$ (the number of POD modes used in the POD reduced order model).
Let $X := H_0^1(\Omega)$, where $\Omega$ is the computational domain.
Let $u(\cdot,t)\in X, t\in[0,T]$ be the weak solution of the weak formulation of the heat equation:
\begin{eqnarray}
(u_t , v)
+ \nu \, (\nabla u , \nabla v)
= (f , v) \,
\quad \forall v \in X.
\label{eqn:heat_weak}
\end{eqnarray}
Given the time instances $t_1, \ldots, t_N \in [0,T]$, we consider the following two ensembles of snapshots:
\begin{eqnarray}
	V^{no\_DQ} 
	&:=& \mbox{span}\left\{ u(\cdot,t_0), \ldots, u(\cdot,t_N) \right\}, 
	\label{snapshots_nodq} \\
	V^{DQ} 
	&:=& \mbox{span}\left\{ u(\cdot , t_0), \ldots, u(\cdot , t_N), 
								\overline{\partial}  u(\cdot , t_1), \ldots, \overline{\partial}  u(\cdot , t_N) \right\},
	\label{snapshots_dq}
\end{eqnarray}
where 
$\overline{\partial} u(t_n) := \frac{u(t_n) - u(t_{n-1})}{\Delta t} , \ n = 1, \ldots, N$
are the time {\it difference quotients (DQs)}.
The two ensembles of snapshots correspond to the two cases investigated in this paper:
(i) with the DQs not included in the snapshots (i.e., $V^{no\_DQ}$); and 
(ii) with the DQs included in the snapshots (i.e., $V^{DQ}$).
As pointed out in Remark 1 in~\cite{KV01}, the ensemble of snapshots $V^{no\_DQ}$ and $V^{DQ}$ yield {\it different} POD bases.
This is clearly illustrated by Figures~\ref{fig:pod_basis_heat}-\ref{fig:pod_basis_burgers}, which display POD basis functions for the heat equation and the Burgers equation, respectively (see Section~\ref{sec:numerical} for details regarding the numerical simulations).

To simplify the presentation, we denote both sets of snapshots (i.e., $V^{no\_DQ}$ and $V^{DQ}$) by 
\begin{equation*}
V = \mbox{span}\left\{ s_1, s_2, \ldots, s_M \right\}, 
\end{equation*}
where $M=N+1$ when $V^{no\_DQ}$ is considered and $M=2N+1$ when $V^{DQ}$ is considered. We use the specific notation (i.e., $V^{no\_DQ}$ or $V^{DQ}$) only when this is necessary.
Let dim $V = d$. 
The POD seeks a low-dimensional basis $\{ \varphi_1, \ldots, \varphi_r \}$,
with $r \leq d$, which optimally approximates the input collection.
Specifically, the POD basis satisfies
\begin{eqnarray}
  \min \frac{1}{M} \sum_{i=1}^M \left\| s_i - 
  \sum_{j=1}^r \bigl( s_i \, , \, \varphi_j(\cdot) \bigr)_{L^2} \, \varphi_j(\cdot) 
  \right\|_{L^2}^2 \, ,
\label{pod_min}
\end{eqnarray}
subject to the conditions that $(\varphi_i,\varphi_j)_{L^2} = \delta_{ij}, \ 1 \leq i, j \leq r$.

In order to solve \eqref{pod_min}, we consider the eigenvalue problem
\begin{eqnarray}
K \, v = \lambda \, v \, ,
\label{pod_eigenvalue}
\end{eqnarray}
where $K \in \R^{M \times M}$, with 
$\displaystyle K_{ij} = \frac{1}{M} \, ( s_j, s_i )_{L^2} \,$, 
is the snapshot correlation matrix,
$\lambda_1 \geq \lambda_2 \geq \ldots \geq \lambda_d >0$ are the positive eigenvalues, and
$v_k, \, k = 1, \ldots, d,$ are the associated eigenvectors.
It can then be shown (see, e.g.,  \cite{HLB96,KV99}), that the solution of \eqref{pod_min} is given by
\begin{eqnarray}
\varphi_{k}(\cdot) = \frac{1}{\sqrt{\lambda_k}} \, \sum_{j=1}^{M} (v_k)_j \, s_{j},
\quad 1 \leq k \leq r,
\label{pod_basis_formula}
\end{eqnarray}
where $(v_k)_j$ is the $j$-th component of the eigenvector $v_k$.

\begin{definition}
The term 
\begin{eqnarray}
\eta^{interp}(x,t)
:= u(x,t) 
- \sum_{j=1}^r \bigl( u(x,t) \, , \, \varphi_j(x) \bigr)_{L^2} \, \varphi_j(x)
\label{eqn:definition_interpolation_error}
\end{eqnarray}
will be denoted as the {\it POD interpolation error}.
\end{definition}
\begin{remark}
\label{remark:pod_norm}
Note that the $H^1$-norm can also be used to generate the POD basis \cite{KV01,singler2012new}.
In this case, $\| \eta^{interp} \|_{L^2} \sim \| \nabla \eta^{interp} \|_{L^2}$. 
Thus, the two cases considered in this paper (i.e., $V = V^{no\_DQ}$ and $V = V^{DQ}$) yield error estimates that have the same convergence rates with respect to $r$.
\end{remark}

It can also be shown~\cite{KV01} that the following {\it POD approximation property} holds:
\begin{eqnarray}
&& \frac{1}{N+1} \sum_{i=0}^N \left\| u(\cdot,t_i) - 
  \sum_{j=1}^r \bigl( u(\cdot,t_i) \, , \, \varphi_j(\cdot) \bigr)_{L^2} \, \varphi_j(\cdot) 
  \right\|_{L^2}^2
= \sum_{j=r+1}^{d} \lambda_j 
\qquad \text{if } \ V = V^{no\_DQ} \, , 
\label{pod_error_formula_nodq} \\
&& \frac{1}{2 N + 1} \, \sum_{i=0}^N 
\left\| 
u(\cdot,t_i) - \sum_{j=1}^r \bigl( u(\cdot,t_i) \, , \, \varphi_j(\cdot) \bigr)_{L^2} \, \varphi_j(\cdot) 
\right\|_{L^2}^2  
\label{pod_error_formula_dq} \\
&& + \  \frac{1}{2 N + 1} \, \sum_{i=1}^N 
\left\| 
\overline{\partial}  u(\cdot,t_i) - \sum_{j=1}^r \bigl( \overline{\partial}  u(\cdot,t_i) \, , \, \varphi_j(\cdot) \bigr)_{L^2} \, \varphi_j(\cdot) 
\right\|_{L^2}^2
= \sum_{j=r+1}^{d} \lambda_j
\quad \text{if } \ V = V^{DQ}  \, . \quad
\nonumber
\end{eqnarray}
The approximation property~\eqref{pod_error_formula_nodq}-\eqref{pod_error_formula_dq} represents the relationship between the average of the square of the $L^2$-norm of the interpolation error and the sum of the eigenvalues of the POD modes that are not included in the POD reduced order model.

%
In order to be able to prove pointwise in time error estimates in Section~\ref{sec:error}, we also make the following assumption:
\begin{assumption}
We assume that, for $i = 1, \ldots, N$, the interpolation error satisfies the the following estimates:
\begin{eqnarray}
&& \left\| u(\cdot,t_i) - 
  \sum_{j=1}^r \bigl( u(\cdot,t_i) \, , \, \varphi_j(\cdot) \bigr)_{L^2} \, \varphi_j(\cdot) 
  \right\|_{L^2}^2
\leq C \, \sum_{j=r+1}^{d} \lambda_j 
\qquad \text{if } \ V = V^{no\_DQ} \, , 
\label{pod_error_formula_nodq_pointwise} \\
&&  \left\| 
u(\cdot,t_i) - \sum_{j=1}^r \bigl( u(\cdot,t_i) \, , \, \varphi_j(\cdot) \bigr)_{L^2} \, \varphi_j(\cdot) 
\right\|_{L^2}^2  
+ \ \left\| 
\overline{\partial}  u(\cdot,t_i) - \sum_{j=1}^r \bigl( \overline{\partial}  u(\cdot,t_i) \, , \, \varphi_j(\cdot) \bigr)_{L^2} \, \varphi_j(\cdot) 
\right\|_{L^2}^2 \nonumber \\
&& \hspace*{6.2cm} \leq C \,  \sum_{j=r+1}^{d} \lambda_j
\quad \text{if } \ V = V^{DQ}  \, . \quad
\label{pod_error_formula_dq_pointwise}
\end{eqnarray}
\label{assumption_pod_error_formula_pointwise}
\end{assumption}

\begin{remark}
Assumption~\ref{assumption_pod_error_formula_pointwise} is natural.
It simply says that in the sums in \eqref{pod_error_formula_nodq} and \eqref{pod_error_formula_dq} no individual term is much larger than the other terms in this sum.
We also note that Assumption~\ref{assumption_pod_error_formula_pointwise} does not play an essential role in the error analysis in Section~\ref{sec:error}, since we will exclusively consider the continuous in time formulation.

We mention that Assumption~\ref{assumption_pod_error_formula_pointwise} would follow directly from the POD approximation property \eqref{pod_error_formula_nodq}--\eqref{pod_error_formula_dq} if we dropped the $\frac{1}{M}$ factor in the snapshot correlation matrix $K$.
In fact, this approach is used in, e.g., \cite{KV02,volkwein2011model}.
We note, however, that this would most likely increase the magnitudes of the eigenvalues on the RHS of the POD approximation property \eqref{pod_error_formula_nodq}--\eqref{pod_error_formula_dq}. 
\label{remark_pod_error_formula_pointwise}
\end{remark}

In what follows, we will use the notation 
$X^r = \mbox{span}\{ \varphi_1, \varphi_2, \ldots, \varphi_r \} \, .$
To derive the POD reduced-order model for the heat equation~\eqref{eqn:heat_weak}, we employ the Galerkin truncation, which yields the following approximation $u_r \in X^r$ of $u$:
\begin{eqnarray}
u_r(x, t)
:= \sum_{j=1}^{r} a_j(t) \, \varphi_j(x) .
\label{pod_g_truncation}
\end{eqnarray}
Plugging \eqref{pod_g_truncation} into \eqref{eqn:heat_weak} and multiplying by test functions in $X^r \subset X$ yields the {\em POD Galerkin reduced order model (POD-G-ROM)}:
\begin{eqnarray}
(u_{r, t} , v_r)
+ \nu \, (\nabla u_r , \nabla v_r)
= (f , v_r)
\quad \forall \, v_r \in X^r .
\label{eqn:pod_g_rom}
\end{eqnarray}
The main advantage of the POD-G-ROM \eqref{eqn:pod_g_rom} over a straightforward finite element discretization of \eqref{eqn:heat_weak} is clear --
the computational cost of the former is dramatically lower than that of the latter.


\section{Error Estimates}
\label{sec:error}

In this section, we prove estimates for the error $ u - u_r $, where $u$ is the solution of the weak formulation of the heat equation \eqref{eqn:heat_weak} and $u_r$ is the solution of the POD-G-ROM~\eqref{eqn:pod_g_rom}. 
Error estimates for the POD reduced order modeling of general systems were derived in, e.g., \cite{KV01,KV02,veroy2005certified,HPS07,LCNY08,rozza2008reduced,kalashnikova2010stability,drohmann2012reduced,urban2012new,amsallem2013error,sachs2013priori,herkt2013convergence}
In our theoretical analysis, we consider two cases, depending on the type of snapshots used in the derivation of the POD basis: {\bf Case I}: $V = V^{DQ}$ (i.e., with the DQs); and {\bf Case II}: $V = V^{no\_DQ}$ (i.e., without the DQs).
The main goal of this paper is to investigate whether {\bf Case I}, {\bf Case II}, or both {\bf Case I} and {\bf Case II}, yield {\it error estimates that are optimal with respect to $r$}.
The optimality with respect to $r$ is given by the following error estimates: 
\begin{eqnarray}
\| u - u_r \| 
&\leq& C \, \| \eta^{interp} \| \, ,
\label{eqn:optimality_eta} \\
\| \nabla (u - u_r) \| 
&\leq& C \, \| \nabla \eta^{interp}\| \, ,
\label{eqn:optimality_nabla_eta}
\end{eqnarray}
where $\eta^{interp}$ is the POD interpolation error defined in~\eqref{eqn:definition_interpolation_error}.

We emphasize that $\| \eta^{interp} \|$ and $\| \nabla \eta^{interp}\|$ scale differently with respect to $r$:
The scaling of $\| \eta^{interp} \|$ is given by the POD approximation property~\eqref{pod_error_formula_nodq_pointwise}--\eqref{pod_error_formula_dq_pointwise} in Assumption~\ref{assumption_pod_error_formula_pointwise}.
The scaling of $\| \nabla \eta^{interp} \|$ is {\it not} given by the POD approximation property~\eqref{pod_error_formula_nodq_pointwise}--\eqref{pod_error_formula_dq_pointwise} in Assumption~\ref{assumption_pod_error_formula_pointwise}.
To derive such an estimate, we use the fact that the interpolation error lives in a finite dimensional space, i.e., the space spanned by the snapshots.
Using an inverse estimate similar to that presented in Lemma~\ref{lemma_inverse_pod} but for the entire space of snapshots (of dimension $d$), we get the following estimate:
\begin{eqnarray}
\| \nabla \eta^{interp} \|_{L^2}
\leq C_{inv}(d) \, \| \eta^{interp} \|_{L^2} \, ,
\label{eqn:inverse_estimate_interpolation_error}
\end{eqnarray}
where $C_{inv}(d)$ is the constant in the inverse estimate in Lemma~\ref{lemma_inverse_pod}.
Following the discussion in Remark~\ref{remark_pod_inverse_estimate_scalings}, we conclude that the scaling of $\| \nabla \eta^{interp} \|$ is of lower order with respect to $r$ than the scaling of $\| \eta^{interp} \|$.
Thus, if the error analysis yields estimates of the form
\begin{eqnarray}
\| u - u_r \| 
\leq C \, \| \nabla \eta^{interp}\| \, ,
\label{eqn:suboptimal_error_estimate}
\end{eqnarray}
these estimates will be called {\it suboptimal} with respect to $r$.

In this section, we investigate the optimality question from a theoretical point of view, by monitoring the dependency of the error estimates on $r$.
In Section~\ref{sec:numerical}, we investigate the same question from a numerical point of view.

We note that we perform the error analysis only for the {\it semidiscretization} of the POD-G-ROM~\eqref{eqn:pod_g_rom}. 
In fact, in this semidiscretization, we only consider the error component corresponding to the POD truncation.
Of course, in practical numerical simulations, the semidiscretization also has a spatial component (e.g., due to the finite element discretization -- see Section~\ref{sec:numerical}).
Furthermore, when considering the full discretization, the error also has a time discretization component (e.g., due to the time stepping algorithm -- see Section~\ref{sec:numerical}).
All these error components should be included in a rigorous error analysis of the discretization of the POD-G-ROM~\eqref{eqn:pod_g_rom} (see, e.g., \cite{LCNY08,iliescu2012variational,iliescu2013variational}).
For clarity of presentation, however, we only consider the error component corresponding to the POD truncation.
In what follows, we will show that this is sufficient for answering the question asked in the title of this paper.

We start by introducing some notation and we list several results that will  be used throughout this section.


We note that, since the POD basis is computed in the $L^2$-norm (see~\eqref{pod_min}), the POD mass matrix $M_{r} \in \R^{r \times r}$ with $M_{i j} = (\varphi_j , \varphi_i)$ is the identity matrix.
Thus, the POD inverse estimate that was proven in \cite{iliescu2012variational} (see also Lemma 2 and Remark 2 in~\cite{KV01})  becomes:
\begin{lemma}[POD Inverse Estimate]
Let $S_{r} \in \R^{r \times r}$ with $S_{i j} = (\nabla \varphi_j , \nabla \varphi_i)$ be the POD stiffness matrix and let $\| \cdot \|_2$  denote the matrix 2-norm.
Then, for all $v_r \in X^r$, the following POD inverse estimate holds:
\begin{eqnarray}
\| \nabla v_r \|_{L^2}  
\leq C_{inv}(r) \, \| v_r \|_{L^2}  \, ,
\label{lemma_inverse_pod_1}
\end{eqnarray}
where $C_{inv}(r) := \sqrt{\| S_r \|_2}$.
\label{lemma_inverse_pod}
\end{lemma}



\begin{remark}[POD Inverse Estimate Scalings]
\label{remark_pod_inverse_estimate_scalings}
Since the $r$ dependency of the error estimates presented in this section will be carefully monitored, we try to get some insight into the scalings of the constant $C_{inv}(r)$ in \eqref{lemma_inverse_pod_1}, i.e., the scalings of $\| S_r \|_2$ with respect to $r$, the dimension of the POD basis.
We note that, since the POD basis significantly varies from test case to test case, it would be difficult to derive a general scaling of $\| S_r \|_2$. 
We emphasize, however, that when the underlying system is {\it homogeneous} (i.e., invariant to spatial translations), the POD basis is identical to the Fourier basis (see, e.g., \cite{HLB96}).
In this case, it is easy to derive the scalings of $\| S_r \|_2$ with respect to $r$.
Without loss of generality, we assume that the computational domain is $[0,1] \subset \R^{1}$ and that the boundary conditions are homogeneous Dirichlet. 
In this case, the Fourier basis functions are given by the following formula: $\varphi_j(x) = \sin(j \, \pi \, x)$.
It is a simple calculation to show that the matrix $S_r$ is diagonal and that its diagonal entries are given by the following formula:
\begin{eqnarray}
S_{j j}
= \int_{0}^{1} ( j \, \pi)^2 \, \sin^2(j \, \pi \, x) \, dx
= \frac{1}{2} \, ( j \, \pi)^2 \, .
\label{eqn:pod_stiffness_matrix_entries_1d}
\end{eqnarray}
It is easy to see that, in the $n$-dimensional case, formula~\eqref{eqn:pod_stiffness_matrix_entries_1d} becomes
\begin{eqnarray}
S_{j j}
= \int_{\Omega} ( j \, \pi)^{2 \, n} \, \sin^2(j \, \pi \, x) \, dx
= \frac{1}{2} \, ( j \, \pi)^{2 \, n} \, .
\label{eqn:pod_stiffness_matrix_entries_nd}
\end{eqnarray}
Since the POD stiffness matrix $S_r$ is symmetric, its matrix 2-norm is given by $\| S_r \|_2 = \lambda_{max}$, where $\lambda_{max}$ is the largest eigenvalue of $S_r$.
Thus, in the $n$-dimensional case, we have
\begin{eqnarray}
\| S_r \|_2 
= \frac{1}{2} \, ( r \, \pi)^{2 \, n}
= \mathcal{O}(r^{2 \, n}) \, .
\label{eqn:pod_stiffness_matrix_scaling_nd}
\end{eqnarray}
Thus, we conclude that, when the underlying system is homogeneous, the 2-norm of the POD stiffness matrix $S_r$ scales as $\mathcal{O}(r^{2 \, n})$, where $n$ is the spatial dimension.

As mentioned at the beginning of the remark, for general (non-homogeneous) systems is would be hard to derive theoretical scalings.
The numerical tests in Section~\ref{sec:numerical}, however, seem to confirm the theoretical scaling in \eqref{eqn:pod_stiffness_matrix_scaling_nd}.

For the heat equation and the Burgers equation (see Section~\ref{sec:numerical} for details regarding the numerical simulations), we monitor the scaling of $\|S_r\|_2$ with respect to $r$ for the POD-G-ROM~\eqref{eqn:pod_g_rom}.
We consider two cases: when the DQs are used in the generation of the POD basis (the corresponding results are denoted by DQ), and when the DQs are not used in the generation of the POD basis (the corresponding results are denoted by no-DQ).
The scalings are plotted in Figure \ref{fig:test1_Sr}. 
It is seen that both no-DQ and DQ approaches yields scalings that are similar to the theoretical scaling~\eqref{eqn:pod_stiffness_matrix_scaling_nd} predicted for the homogeneous flow fields (i.e., when the POD basis reduces to the Fourier basis).
The only exception seems to be for the Burgers equation in the DQ case.
In all cases, however, the scaling $\| S_r \|_2 = \mathcal{O}(r^{\alpha})$, where $\alpha$ is a positive constant, seems to be valid.

\begin{figure}[htp]
\begin{minipage}[h]{.45\linewidth}
\includegraphics[width=1.1\textwidth]{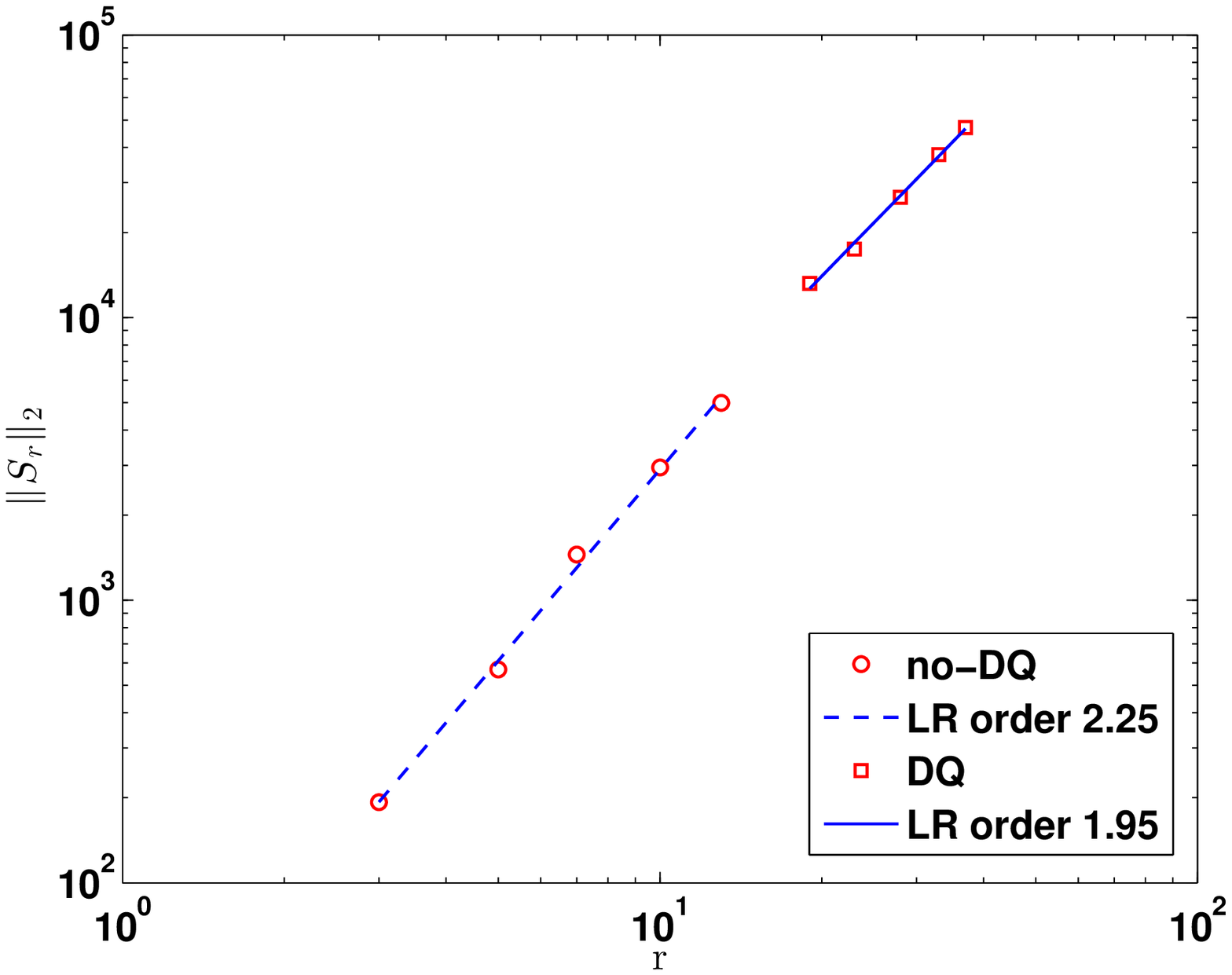}
\end{minipage}
\hspace{.3cm}
\begin{minipage}[h]{.45\linewidth}
\includegraphics[width=1.1\textwidth]{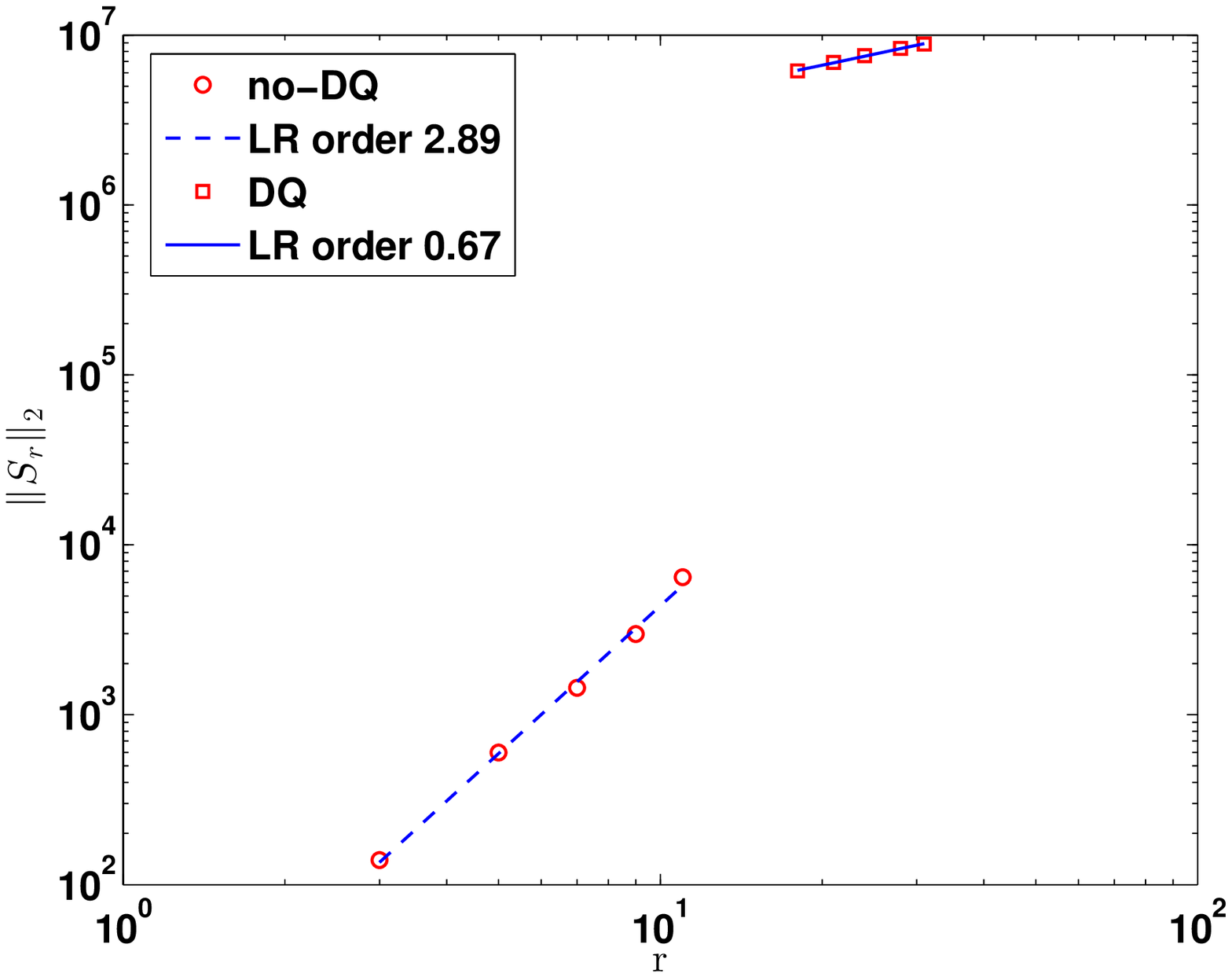}
\end{minipage}
\caption{
	Heat equation (left); Burgers equation (right).
	Plots of the scalings of $\|S_r\|_2$ with respect to $r$ when the DQs are used (denoted by DQ) and when the DQs are not used (denoted by no-DQ).
	}
\label{fig:test1_Sr}
\end{figure}

\end{remark}

\medskip

After these preliminaries, we are ready to derive the error estimates.
The error analysis will proceed along the same lines as the error analysis of the finite element semidiscretization~\cite{GR86,layton2008introduction,thomee2006galerkin}.
The main difference between the finite element and the POD settings is that the finite element approximation property is {\it global} \cite{GR86,layton2008introduction}, whereas the POD approximation property is {\it local}, i.e., it is only valid at the time instances at which the snapshots were taken (see~\eqref{pod_error_formula_nodq}-\eqref{pod_error_formula_dq_pointwise}).
Thus, in order to be able to use the POD approximation property~\eqref{pod_error_formula_nodq}-\eqref{pod_error_formula_dq_pointwise}, in what follows we assume that the final error estimates for the semidiscretization are considered only at the time instances $t_i, \, i =1, \ldots, N$.

We start by considering the error equation:
\begin{eqnarray}
( e_t , v_r )
+ \nu \, ( \nabla e , \nabla v_r )
= 0 
\qquad
\forall \, v_r \in X^r ,
\label{s_error_semi_eq_1}
\end{eqnarray}
where 
$e := u - u_r$ is the error.
The error is split into two parts:
\begin{eqnarray}
e
= u - u_r
= ( u - w_r )
-  ( u_r - w_r ) 
= \eta
- \phi_r ,
\label{s_error_semi_eq_2}
\end{eqnarray}
where 
$w_r$ is an {\it arbitrary} function in $X^r$,
$\eta := u - w_r$, and
$\phi_r := u_r - w_r$.
Using this decomposition in the error equation \eqref{s_error_semi_eq_1}, we get
\begin{eqnarray}
( \phi_{r,t} , v_r )
+ \nu \, ( \nabla \phi_r , \nabla v_r )
= ( \eta_t , v_r )
+ \nu \, ( \nabla \eta , \nabla v_r ) .
\label{s_error_semi_eq_2a}
\end{eqnarray}
The analysis proceeds by using \eqref{s_error_semi_eq_2a} to show that
\begin{eqnarray}
\| \phi_r \|
\leq C \, \| \eta \| .
\label{s_error_semi_eq_3}
\end{eqnarray}
Using the triangle inequality, one then gets
\begin{eqnarray}
\| e \|
\leq \| \eta \|
+ \| \phi_r \|
\leq (1 + C) \, \| \eta \| .
\label{s_error_semi_eq_4}
\end{eqnarray}
Since $w_r$ was chosen arbitrarily and since the LHS does not depend on $w_r$, we can take the infimum over $w_r$ in \eqref{s_error_semi_eq_4} and get the following error estimate:
\begin{eqnarray}
\| e \|
\leq (1 + C) \, \inf_{w_r \in X^r} \| u - w_r \| .
\label{s_error_semi_eq_5}
\end{eqnarray}
At this stage, one invokes the approximability property of the space $X^r$ in \eqref{s_error_semi_eq_5} and concludes that the error estimate~\eqref{s_error_semi_eq_5} is {\it optimal with respect to $r$}.
That is, the error in~\eqref{s_error_semi_eq_5} is, up to a constant, the interpolation error in $X^r$.

There are a lot of variations of the error analysis sketched above, but most of the existing error analyses follow the above path.
The main variations are the following:
(i) the choice of the arbitrary function $w_r \in X^r$;
(ii) the norm $\| \cdot \|$ used above; and
(iii) the choice of the test function used in the error equation to get \eqref{s_error_semi_eq_3} \cite{thomee2006galerkin}.

In the remainder of this section, we investigate whether error estimates that are {\it optimal with respect to $r$} can be obtained with or without including the DQs in the set of snapshots.
To this end, in Section~\ref{sec:case_I} we consider the case in which the DQs are included in the set of snapshots (i.e., $V = V^{DQ}$).
Then, in Section~\ref{sec:case_II} we consider the case in which the DQs are not included in the set of snapshots (i.e., $V = V^{no\_DQ}$).

\subsection{Case I ($V = V^{DQ}$)}
\label{sec:case_I}

The standard approach used to prove error estimates in this case is to use the Ritz projection~\cite{KV01,KV02,iliescu2012variational,iliescu2013variational}.
We note that this is the standard approach used in the finite element context~\cite{GR86,layton2008introduction}.
As pointed out in \cite{thomee2006galerkin}, the Ritz projection was first used by Wheeler in \cite{wheeler1973priori} to obtain {\it optimal} error estimates for the finite element discretization of parabolic problems.

We start by choosing $w_r := R_r(u)$ in~\eqref{s_error_semi_eq_2}, where $R_r(u)$ is the {\it Ritz projection} of $u$, given by:
\begin{eqnarray}
\bigl( \nabla ( u - R_r(u) ) , \nabla v_r \bigr)
= 0
\qquad \forall \, v_r \in X^r .
\label{s_error_semi_eq_6}
\end{eqnarray}
To emphasize that we are using the Ritz projection, in the remainder of Section~\ref{sec:case_I} we will use the notation $\eta = u - R_r(u) = \eta^{Ritz}$.

Using \eqref{s_error_semi_eq_6}, \eqref{s_error_semi_eq_2a} becomes
\begin{eqnarray}
( \phi_{r,t} , v_r )
+ \nu \, ( \nabla \phi_r , \nabla v_r )
= ( \eta^{Ritz}_t , v_r )
+ \nu \, \cancelto{0}{( \nabla \eta^{Ritz} , \nabla v_r )} .
\label{s_error_semi_eq_7}
\end{eqnarray}
It is the cancelation of the last term on the RHS of \eqref{s_error_semi_eq_7} that yields {\it optimal} error estimates.
We let $v_r := \phi_r$ in \eqref{s_error_semi_eq_7}, and then we apply the Cauchy-Schwarz inequality to the remaining term on the RHS:
\begin{eqnarray}
\frac{1}{2} \, \frac{d}{dt} \| \phi_r \|^2
+ \nu \, \| \nabla \phi_r \|^2
\leq \| \eta^{Ritz}_t \| \, \| \phi_r \| .
\label{s_error_semi_eq_8}
\end{eqnarray}
We rewrite the first term on the LHS of \eqref{s_error_semi_eq_8} as
\begin{eqnarray}
\frac{1}{2} \, \frac{d}{dt} \| \phi_r \|^2
= \| \phi_r \| \, \frac{d}{dt} \| \phi_r \| .
\label{s_error_semi_eq_9}
\end{eqnarray}
We apply the Poincare-Friedrichs inequality to the second term on the LHS of \eqref{s_error_semi_eq_8}:
\begin{eqnarray}
\nu \, \| \nabla \phi_r \|^2
\geq C \, \nu \, \| \phi_r \|^2 .
\label{s_error_semi_eq_10}
\end{eqnarray}
We note that the Poincare-Friedrichs inequality $C \| v \|^2 \leq \| \nabla v \|^2$ holds for every function $v$ in the {\it continuous} space $H_0^1(\Omega)$, and, in particular, for $\phi_r \in X^r \subset X = H_0^1(\Omega)$ (see equation (3) in \cite{KV01}).
Thus, the constant $C$ in \eqref{s_error_semi_eq_10} does not depend on $r$.
Using \eqref{s_error_semi_eq_9} and \eqref{s_error_semi_eq_10} in \eqref{s_error_semi_eq_8}, we get
\begin{eqnarray}
\frac{1}{2} \, \frac{d}{dt} \| \phi_r \|
+ C \, \nu \, \| \phi_r \|
\leq \| \eta^{Ritz}_t \| .
\label{s_error_semi_eq_11}
\end{eqnarray}
Using Gronwall's lemma in \eqref{s_error_semi_eq_11}, we get for $0 < t \leq T$
\begin{eqnarray}
\| \phi_r (t) \|
\leq e^{-2 \, C \, \nu \, t} \, \| \phi_r(0) \|
+ 2 \, \int_0^t e^{-2 \,C \, \nu \, (t-s)} \, \| \eta^{Ritz}_t(s) \| \, ds .
\label{s_error_semi_eq_12}
\end{eqnarray}
Using~\eqref{s_error_semi_eq_4}, the first term on the RHS of \eqref{s_error_semi_eq_12} can be estimated as follows:
\begin{eqnarray}
\| \phi_r(0) \|
\leq \| e(0) \| 
+ \| \eta^{Ritz}(0) \| .
\label{s_error_semi_eq_13}
\end{eqnarray}
Thus, \eqref{s_error_semi_eq_12} becomes
\begin{eqnarray}
\| \phi_r (t) \|
\leq e^{-2 \, C \, \nu \, t} \, \biggl( \| e(0) \| 
+ \| \eta^{Ritz}(0) \| \biggr)
+ 2 \, \int_0^t e^{-2 \, C \, \nu \, (t-s)} \, \| \eta^{Ritz}_t(s) \| \, ds .
\label{s_error_semi_eq_14}
\end{eqnarray}
Applying the triangle inequality, just as in \eqref{s_error_semi_eq_4}, we get
\begin{eqnarray}
\| e(t) \|
\leq  \| \eta^{Ritz}(t) \|
+ e^{-2 \, C \, \nu \, t} \, \biggl( \| e(0) \| 
+ \| \eta^{Ritz}(0) \| \biggr)
+ 2 \, \int_0^t e^{-2 \, C \, \nu \, (t-s)} \, \| \eta^{Ritz}_t(s) \| \, ds . 
\nonumber \\
\label{s_error_semi_eq_15}
\end{eqnarray}
Estimate~\eqref{s_error_semi_eq_15} shows that, as long as the Ritz projection~\eqref{s_error_semi_eq_6} yields estimates for $\| \eta^{Ritz}(t) \|$ (including for $t=0$) and $\| \eta^{Ritz}_t(t) \|$ that are {\it optimal} with respect to $r$, the estimates for the POD error $e$ are also optimal with respect to $r$.
\begin{remark}[POD Ritz Projection]
\label{remark_ritz}
In the finite element context, both $\| \eta^{Ritz} \|$ and $\| \eta^{Ritz}_t \|$ are optimal with respect to the mesh size $h$.
In the POD context, however, this is not that clear.
To the best of the authors' knowledge, the state-of-the-art regarding the Ritz projection in a POD context is given in the pioneering paper of Kunisch and Volkwein~\cite{KV01}.
Since the Ritz projection plays such an important role in this paper, we summarize below the results in~\cite{KV01}. 


The main result in \cite{KV01} regarding the Ritz projection is Lemma 3 (see also (10) and (11) in \cite{KV01}), which, in our notation, states the following:
\begin{eqnarray}
\| \nabla \eta^{Ritz} \|^ 2 
\leq \| S_d \|_2 \, \sum \limits_{j=r+1}^{d} \lambda_j \, .
\label{eqn:kv01_lemma3_2}
\end{eqnarray}
For clarity of presentation, we have not included in \eqref{eqn:kv01_lemma3_2} the constants that do not depend on $r$.
Thus, the following relationship between the Ritz projection error and the POD interpolation error holds:
\begin{eqnarray}
\| \nabla \eta^{Ritz} \|
\leq C \, \sqrt{\| S_d \|_2} \, \| \eta^{interp} \| \, .
\label{eqn:relationship_ritz_interpolation_nabla_eta}
\end{eqnarray}
The scaling in \eqref{eqn:relationship_ritz_interpolation_nabla_eta} suggests that the Ritz projection yields optimal error estimates with respect to $r$ in the $H^1$-seminorm (see \eqref{eqn:inverse_estimate_interpolation_error}).

We emphasize that Lemma 3 in \cite{KV01} does not include any bounds for the $L^2$-norm of $\eta^{Ritz}$, i.e., for $\| \eta^{Ritz} \|$.
This is in clear contrast with the finite element context, in which $\| \eta^{Ritz} \|$ is estimated by the usual duality argument (the Aubin-Nitsche ``trick," see, e.g., \cite{thomee2006galerkin}).
Using a duality argument, however, is challenging in the POD context, since any auxiliary dual problem would not necessarily inherit the POD approximation property~\eqref{pod_error_formula_nodq}--\eqref{pod_error_formula_dq_pointwise}.
To the best of the authors' knowledge, such a duality argument has never been used in a POD context.

We emphasize that not being able to use a duality argument in the Ritz projection to get error estimates that are optimal with respect to $r$ has significant consequences in the error analysis.
Indeed, in the proof of Theorem 7 in \cite{KV01} (the error estimate for the backward Euler time discretization), to estimate the $\| \eta^{Ritz} \|$ error component in equations (27a) and (27b), the authors use the $\| \nabla \eta^{Ritz} \|$ estimate given in Lemma 3 and the Poincare-Friedrichs inequality given in equation (3).
Since the Poincare-Friedrichs constant does not depend on $r$, we conclude that $\| \eta^{Ritz}\|$ and $\| \nabla \eta^{Ritz} \|$ have the same order. 
This, in turn, suggests that $\| \eta^{Ritz} \|$ is suboptimal with respect to $r$.
We note that the same approach (i.e., Lemma 4 and the Poincare-Friedrichs inequality) is used in~\cite{KV01} to estimate the DQ approximation of $\| \eta^{Ritz}_t \|$ (see the two inequalities above equation (29a)). 

Thus, the analysis in \cite{KV01} suggests that the estimates for $\| \eta^{Ritz} \|$ and $\| \eta^{Ritz}_t \|$ are suboptimal with respect to $r$.
The numerical tests in Section~\ref{sec:numerical} of this report, however, suggest that these estimates are actually optimal, just as in the finite element case.
In the analysis that follows, we will use the insight from the numerical results in Section~\ref{sec:numerical} and 
make the following assumption:
\begin{assumption}
We assume that the POD Ritz projection error $\eta^{Ritz}$ satisfies optimal error estimates with respect to $r$ in the $L^2$-norm:\begin{eqnarray}
\| \eta^{Ritz} \|
&\leq& C \, \| \eta^{interp} \| \, ,
\label{eqn:relationship_ritz_interpolation_eta} \\
\| \eta^{Ritz}_t \|
&\leq& C \, \| \eta^{interp} \| \, .
\label{eqn:relationship_ritz_interpolation_eta_t}
\end{eqnarray}
\label{assumption_ritz}
\end{assumption}
\end{remark}



\subsection{Case II ($V = V^{no\_DQ}$)}
\label{sec:case_II}

This approach was used in \cite{chapelle2012galerkin,singler2012new}.
The motivation for this approach is the following:
In Case I ($V = V^{DQ}$), the first term on the RHS of \eqref{s_error_semi_eq_7}, $( \eta_t , v_r)$, yields a term $\| \eta_t \|$ that stays in all the subsequent inequalities, including the final error estimate \eqref{s_error_semi_eq_15}. 
Chapelle et al. proposed in \cite{chapelle2012galerkin} a different approach that eliminated the $( \eta_t , v_r)$ term in \eqref{s_error_semi_eq_7}.
Their approach was straightforward: Instead of using the Ritz projection (as in Case I), they used the $L^2$ projection.
That is, they chose $w_r := P_r(u)$, where $P_r(u)$ is the {\it $L^2$ projection} of $u$, given by
\begin{eqnarray}
\bigl( u - P_r(u) , v_r \bigr)
= 0
\qquad \forall \, v_r \in X^r .
\label{s_error_semi_eq_16}
\end{eqnarray}
To emphasize that we are using the $L^2$ projection, in the remainder of Section~\ref{sec:case_II} we will use the notation $\eta = u - P_r(u) = \eta^{L^2}$.

We note that the POD $L^2$ projection error $\eta^{L^2}$ is exactly the POD interpolation error defined in~\eqref{eqn:definition_interpolation_error}:
\begin{eqnarray}
\eta^{L^2}
= \eta^{interp} \, .
\label{s_error_semi_eq_16.5}
\end{eqnarray}

Next, we show how the error analysis in Case I changes with $w_r = P_r(u)$ as in \cite{chapelle2012galerkin} (see also \cite{singler2012new}).

Using \eqref{s_error_semi_eq_16}, \eqref{s_error_semi_eq_2a} becomes
\begin{eqnarray}
( \phi_{r,t} , v_r )
+ \nu \, ( \nabla \phi_r , \nabla v_r )
= \cancelto{0}{( \eta^{L^2}_t , v_r )}
+ \nu \, ( \nabla \eta^{L^2} , \nabla v_r ) .
\label{s_error_semi_eq_17}
\end{eqnarray}

We emphasize that it is the cancelation of the first term on the RHS of \eqref{s_error_semi_eq_17} that yields error estimates that do not require the DQs.
We let $v_r := \phi_r$ in \eqref{s_error_semi_eq_17} and we apply Cauchy-Schwarz inequality to the remaining term on the RHS:
\begin{eqnarray}
\frac{1}{2} \, \frac{d}{dt} \| \phi_r \|^2
+ \nu \, \| \nabla \phi_r \|^2
\leq \nu \, \| \nabla \eta^{L^2} \| \, \| \nabla \phi_r \| .
\label{s_error_semi_eq_18}
\end{eqnarray}
The error analysis can then proceed in several directions.

\subsubsection{Approach II.A}
\label{sec:approach_II.A}

One approach is to use Young's inequality in \eqref{s_error_semi_eq_18} to get
\begin{eqnarray}
\frac{1}{2} \, \frac{d}{dt} \| \phi_r \|^2
+ \nu \, \| \nabla \phi_r \|^2
\leq \frac{\nu}{2} \, \| \nabla \eta^{L^2} \|^2 
+ \frac{\nu}{2} \, \| \nabla \phi_r \|^2 ,
\label{s_error_semi_eq_19}
\end{eqnarray}
which implies
\begin{eqnarray}
\frac{1}{2} \, \frac{d}{dt} \| \phi_r \|^2
+ \frac{\nu}{2} \, \| \nabla \phi_r \|^2
\leq \frac{\nu}{2} \, \| \nabla \eta^{L^2} \|^2 .
\label{s_error_semi_eq_20}
\end{eqnarray}
Noticing that the second term on the LHS of \eqref{s_error_semi_eq_20} is positive, we get 
\begin{eqnarray}
\frac{d}{dt} \| \phi_r \|^2
\leq \nu \, \| \nabla \eta^{L^2} \|^2 .
\label{s_error_semi_eq_21}
\end{eqnarray}
Using~\eqref{s_error_semi_eq_16.5} and~\eqref{eqn:suboptimal_error_estimate}, we conclude that Approach II.A will yield error estimates that are suboptimal with respect to $r$.
We also note in passing that estimate~\eqref{s_error_semi_eq_21} suggests that Approach II.A will yield error estimates that are suboptimal with respect to $h$ as well.

%
%

\subsubsection{Approach II.B}
\label{sec:approach_II.B}

The other way of continuing from \eqref{s_error_semi_eq_18} is to apply the POD inverse estimate~\eqref{lemma_inverse_pod_1}:
\begin{eqnarray}
\| \nabla \phi_r \| 
\leq C_{inv}(r) \, \| \phi_r \| ,
\label{s_error_semi_eq_22}
\end{eqnarray}
where 
$C_{inv}(r) :=  \sqrt{\| S_r \|_2 }$.
Using~\eqref{s_error_semi_eq_22} in \eqref{s_error_semi_eq_18} yields
\begin{eqnarray}
\frac{1}{2} \, \frac{d}{dt} \| \phi_r \|^2
+ \nu \, \| \nabla \phi_r \|^2
\leq C_{inv}(r) \, \nu \, \| \nabla \eta^{L^2} \| \, \| \phi_r \| .
\label{s_error_semi_eq_23}
\end{eqnarray}
Dropping $\nu \, \| \nabla \phi_r \|^2$ in \eqref{s_error_semi_eq_23} and simplifying the resulting inequality by $\| \phi_r \|$ (as we did in \eqref{s_error_semi_eq_12}), we get
\begin{eqnarray}
\frac{1}{2} \, \frac{d}{dt} \| \phi_r \|
\leq C_{inv}(r) \, \nu \, \| \nabla \eta^{L^2} \| .
\label{s_error_semi_eq_24}
\end{eqnarray}
Comparing estimate \eqref{s_error_semi_eq_24} with estimate \eqref{s_error_semi_eq_21} in Approach II.A, we note that both estimates have $\| \nabla \eta \|$ on the RHS.
In addition, estimate \eqref{s_error_semi_eq_24} has $C_{inv}(r)$ on the RHS, which increases the suboptimality with respect to $r$ (see Remark~\ref{remark_pod_inverse_estimate_scalings}).
Thus, estimate~\eqref{s_error_semi_eq_24} suggests that Approach II.B yields estimates that are suboptimal with respect to $r$ (and $h$), just as Approach II.A.


\subsubsection{Approach II.C}
\label{sec:approach_II.C}

Since both Approach II.A and Approach II.B yield error estimates that are suboptimal with respect to $r$ in the $L^2$-norm, one can try instead to prove optimal error estimates in the $H^1$-seminorm.
To this end, we use the approach in \cite{thomee2006galerkin} and, instead of choosing $v_r := \phi_r$ in \eqref{s_error_semi_eq_17}, we choose $v_r := \phi_{r,t}$: 
\begin{eqnarray}
\| \phi_{r,t} \|^2
+ \frac{\nu}{2} \, \frac{d}{dt} \| \nabla \phi_r \|^2
\leq \nu \,\| \nabla \eta^{L^2} \| \, \| \nabla \phi_{r,t} \| .
\label{s_error_semi_eq_25}
\end{eqnarray}
Applying Young's inequality and the POD inverse estimate~\eqref{s_error_semi_eq_22} in \eqref{s_error_semi_eq_25}, we get
\begin{eqnarray}
\| \phi_{r,t} \|^2
+ \frac{\nu}{2} \, \frac{d}{dt} \| \nabla \phi_r \|^2
&\leq& \nu \, \| \nabla \eta^{L^2} \| \, \| \nabla \phi_{r,t} \| 
\nonumber \\
&\leq& \frac{\nu^2}{2} \, C_{inv}(r)^2 \, \| \nabla \eta^{L^2} \|^2
+ \frac{1}{2 \, C_{inv}(r)^2} \, \| \nabla \phi_{r,t} \|^2
\nonumber \\
&\stackrel{\eqref{s_error_semi_eq_22}}{\leq}& \frac{\nu^2}{2} \, C_{inv}(r)^2 \, \| \nabla \eta^{L^2} \|^2
+ \frac{1}{2} \, \| \phi_{r,t} \|^2 ,
\label{s_error_semi_eq_26}
\end{eqnarray}
which implies
\begin{eqnarray}
\frac{d}{dt} \| \nabla \phi_r \|^2
\leq  \nu \, C_{inv}(r)^2 \, \| \nabla \eta^{L^2} \|^2 .
\label{s_error_semi_eq_27}
\end{eqnarray}
In contrast with estimate~\eqref{s_error_semi_eq_21} in Approach II.A and estimate~\eqref{s_error_semi_eq_24} in Approach II.B, estimate~\eqref{s_error_semi_eq_27} seems to yield error estimates that are optimal with respect to $r$.
As in estimates~\eqref{s_error_semi_eq_21} and \eqref{s_error_semi_eq_27}, estimate~\eqref{s_error_semi_eq_27} contains the term $\| \nabla \eta^{L^2} \|$ on the RHS.
We note, however, that this term does not cause any problems, since now we are considering the $H^1$-seminorm of the error.
The factor $C_{inv}(r)^2$ in~\eqref{s_error_semi_eq_27}, however, increases the suboptimality with respect to $r$ (see Remark~\ref{remark_pod_inverse_estimate_scalings}).
Thus, estimate~\eqref{s_error_semi_eq_27} suggests that Approach II.C yields estimates that are suboptimal with respect to $r$, just as Approaches II.A and II.B.

Since for {\bf Case I ($V = V^{DQ}$)} in Section~\ref{sec:case_I}, we did not prove error estimates in the $H^1$-norm, for a fair comparison with Approach II.C, we prove these error estimates below.
To this end, we let $v_r := \phi_{r,t}$ in \eqref{s_error_semi_eq_7}:
\begin{eqnarray}
( \phi_{r,t} , \phi_{r,t} )
+ \nu \, ( \nabla \phi_r , \nabla \phi_{r,t} )
= ( \eta^{Ritz}_t , \phi_{r,t} ) .
\label{s_error_semi_eq_28}
\end{eqnarray}
Applying Young's inequality on the RHS of \eqref{s_error_semi_eq_28}, we get
\begin{eqnarray}
\| \phi_{r,t} \|^2
+ \frac{\nu}{2} \, \frac{d}{dt} \| \nabla \phi_r \|^2
&\leq& \| \eta^{Ritz}_t \| \, \| \phi_{r,t} \| 
\nonumber \\
&\leq& \frac{1}{2} \, \| \eta^{Ritz}_t \|^2
+ \frac{1}{2} \, \| \phi_{r,t} \|^2 ,
\label{s_error_semi_eq_29}
\end{eqnarray}
which implies
\begin{eqnarray}
\frac{d}{dt} \| \nabla \phi_r \|^2
\leq \frac{1}{\nu} \, \| \eta^{Ritz}_t \|^2 .
\label{s_error_semi_eq_30}
\end{eqnarray}
Comparing \eqref{s_error_semi_eq_30} with \eqref{s_error_semi_eq_27} in Approach II.C, we note that the latter does not contain the  factor $C_{inv}(r)^2$.
Thus, {\bf Case I ($V = V^{DQ}$)} in Section~\ref{sec:case_I} yields optimal error estimates with respect to $r$, as opposed to Approach II.C.

We also note that, at first glance, estimate~\eqref{s_error_semi_eq_30} suggests that one can get superconvergence in the $H^1$-seminorm.
As mentioned in~\cite{thomee2006galerkin}, however, when applying the triangle inequality one obtains the expected convergence rate:
\begin{eqnarray}
\| \nabla e \|
\leq \| \nabla \eta^{Ritz} \|
+ \| \nabla \phi_r \| .
\label{s_error_semi_eq_31}
\end{eqnarray}

\subsubsection{Approach II.D}
\label{sec:approach_II.D}

Approaches II.A, II.B, and II.C suggest that the pointwise in time (i.e., in the $C^0(0,T; L^2(\Omega))$-norm and  the $C^0(0,T; H^1(\Omega))$-norm) error estimates are suboptimal with respect to $r$.
Thus, we derive error estimates in the solution norm (i.e., in the $L^2(0,T; H^1(\Omega))$-norm).
Integrating~\eqref{s_error_semi_eq_20} from $0$ to $T$, we get
\begin{eqnarray}
\| \phi_r(T) \|^2
+ \frac{\nu}{2} \, \int_{0}^{T} \| \nabla \phi_r(s) \|^2 \, ds
\leq  \| \phi_r(0) \|^2
+ \frac{\nu}{2} \, \int_{0}^{T} \| \nabla \eta^{L^2}(s) \|^2 \, ds .
\label{s_error_semi_eq_23b}
\end{eqnarray}
Using \eqref{s_error_semi_eq_16.5}, \eqref{eqn:optimality_eta}, and \eqref{eqn:optimality_nabla_eta}, we conclude that estimate in \eqref{s_error_semi_eq_23b} is optimal with respect to $r$.
We note that Proposition 3.3 in \cite{chapelle2012galerkin} yields a similar estimate.

%
%
%

\bigskip

%
%
Case II (i.e., $V = V^{no\_DQ}$) yields the following general conclusions when the $L^2$ projection is used instead of the standard Ritz projection:
If the error is computed pointwise in time (i.e., in the $C^0(0,T; L^2(\Omega))$-norm and the $C^0(0,T; H^1(\Omega))$-norm), then the error estimates are suboptimal with respect to $r$.
This is the main message of Approaches II.A, II.B, and II.C.
If, however, the error is computed in the solution norm (i.e., in the $L^2(0,T; H^1(\Omega))$-norm), then the error estimates are optimal with respect to $r$.
This is the main message of Approach II.D.

\bigskip

\section{Numerical Results}
\label{sec:numerical}

The main goal of this section is to numerically investigate the rates of convergence with respect to $r$ of the POD-G-ROM~\eqref{eqn:pod_g_rom} in the two cases considered in Section~\ref{sec:error}: {\bf Case I} ($V = V^{DQ}$) and {\bf Case II} ($V = V^{no\_DQ}$).
Although the error analysis in Section~\ref{sec:error} has been centered around the (linear) heat equation, in this section we consider both the heat equation (Section~\ref{sec:heat}) and the nonlinear Burgers equation (Section~\ref{sec:burgers}).

To measure the errors in the two cases (i.e., $V = V^{DQ}$ and $V = V^{no\_DQ}$), the same norms as those used in Section~\ref{sec:error} are used in this section.
Denoting the error at time $t_j$ by $e_j := u^r_h(\cdot, t_j) - u(\cdot, t_j)$, the following norms are considered:
the error in the $C^0(0,T; L^2(\Omega))$-norm, approximated by
$\cE_{C^0(L^2)} = \max\limits_{0\leq j \leq N} \|e_j\|_{{L^2}(\Om)}$;
the error in the $C^0(0,T; H^1(\Omega))$-norm, approximated by
$\cE_{C^0(H^1)} = \max\limits_{0\leq j \leq N} \|e_j\|_{{H^1}(\Om)}$; 
and the error in the $L^2(0,T; H^1(\Omega))$-norm, approximated by
$\cE_{L^2(H^1)}=\sqrt{\frac{1}{N+1} \sum\limits_{0\leq j \leq N} \| e_j \|_{{H^1}(\Om)}^2}$.
For clarity, we also use the following notation: $\Lambda_r= \sqrt{\sum\limits_{j=r+1}^d \lambda_j}$. 

As mentioned at the beginning of Section~\ref{sec:error}, the POD-G-ROM~\eqref{eqn:pod_g_rom} error estimates are optimal if the following statements hold:

(i) The $L^2$-norm of the error scales as $L^2$-norm of the POD interpolation error.
Using \eqref{eqn:optimality_eta} and \eqref{pod_error_formula_nodq}--\eqref{pod_error_formula_dq_pointwise}, this statement is equivalent to 
\begin{eqnarray}
\cE_{C^0(L^2)}
= \mathcal{O}\left(\sqrt{\sum\limits_{j=r+1}^d \lambda_j} \ \right) \, ,
\label{eqn:numerical_1}
\end{eqnarray}

(ii) The $H^1$-norm of the error scales as $H^1$-norm of the POD interpolation error.
Using \eqref{eqn:optimality_nabla_eta}, \eqref{eqn:inverse_estimate_interpolation_error}, and \eqref{pod_error_formula_nodq}--\eqref{pod_error_formula_dq_pointwise}, this statement is equivalent to 
\begin{eqnarray}
\cE_{C^0(H^1)}
\sim \cE_{L^2(H^1)}
= \mathcal{O}\left(\sqrt{ \| S_d \|_2 \, \sum\limits_{j=r+1}^d \lambda_j}\right) \, ,
\label{eqn:numerical_2}
\end{eqnarray}

Based on the error analysis in Section~\ref{sec:error}, we expect the following convergence rates with respect to $r$:
\begin{table}[h]
\centering
\tabcaption{Theoretical convergence rates for the $no\_DQ$ and the $DQ$ cases.}
\label{tab:numerical_1}
\begin{tabular}{|c|c|c|}
\hline
& & \\[-0.2cm]
& $no\_DQ$ & $DQ$ \\[0.2cm]
 \hline
 & & \\[-0.2cm]
 $\cE_{C^0(L^2)}$ & suboptimal & optimal \\
  & \eqref{s_error_semi_eq_21}; \eqref{s_error_semi_eq_24} & \eqref{s_error_semi_eq_15} \\[0.2cm]
 \hline
 & & \\[-0.2cm]
 $\cE_{C^0(H^1)}$ & suboptimal & optimal \\
  & \eqref{s_error_semi_eq_27} & \eqref{s_error_semi_eq_31} \\[0.2cm]
 \hline
 & & \\[-0.2cm]
 $\cE_{L^2(H^1)}$ & optimal & optimal \\
  & \eqref{s_error_semi_eq_23b} &  \\[0.2cm]
 \hline
 \end{tabular}
 \end{table}

\subsection{Heat Equation}
\label{sec:heat}

\begin{figure}[h]
\centering
\begin{minipage}[h]{.6\linewidth}
\includegraphics[width=1\textwidth]{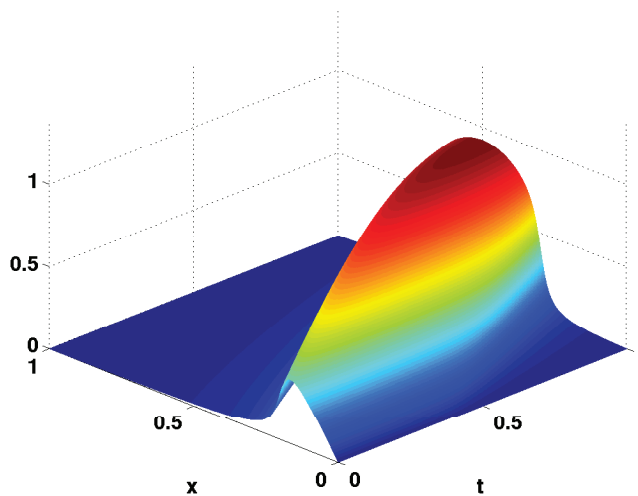}
\end{minipage}
\caption{
	Heat equation.
	Fine resolution finite element solution used to generate the snapshots.
}
\label{fig:test1dns}
\end{figure}

We consider the one-dimensional heat equation~\eqref{eqn:heat_weak} with a known exact solution that represents the propagation in time of a steep front: 
\begin{eqnarray}
u(x,t) = \sin(\pi x)\left[\frac{1}{\pi}\arctan\left(\frac{c}{25}-c\left(x-\frac{t}{2}\right)^2\right) + \frac{1}{2}\right],
\label{eqn:exact_solution_heat}
\end{eqnarray}
where $x\in [0, 1]$ and $t\in [0, 1]$.
The constant $c$ in~\eqref{eqn:exact_solution_heat} controls the steepness of the front.
In all the numerical tests in this section, we used the value $c = 100$.
The value of the diffusion coefficient used in the heat equation~\eqref{eqn:heat_weak} is $\nu = 10^{-2}$. 
Piecewise linear finite elements are used to generate snapshots for the POD-G-ROM~\eqref{eqn:pod_g_rom}. 
A mesh size $h=1/1024$ and the Crank-Nicolson scheme with a time step $\Delta t = 10^{-3}$ are employed for the spatial and temporal discretizations. 
The time evolution of the finite element solution is shown in Figure \ref{fig:test1dns}. 
In total, 1001 snapshots are collected and used for generating POD basis functions in the $L^2$ space. 
The same numerical solver as that used in the finite element approximation is utilized in the POD-G-ROM. 

\begin{figure}[h]
\begin{minipage}[h]{.45\linewidth}
\includegraphics[width=1.1\textwidth]{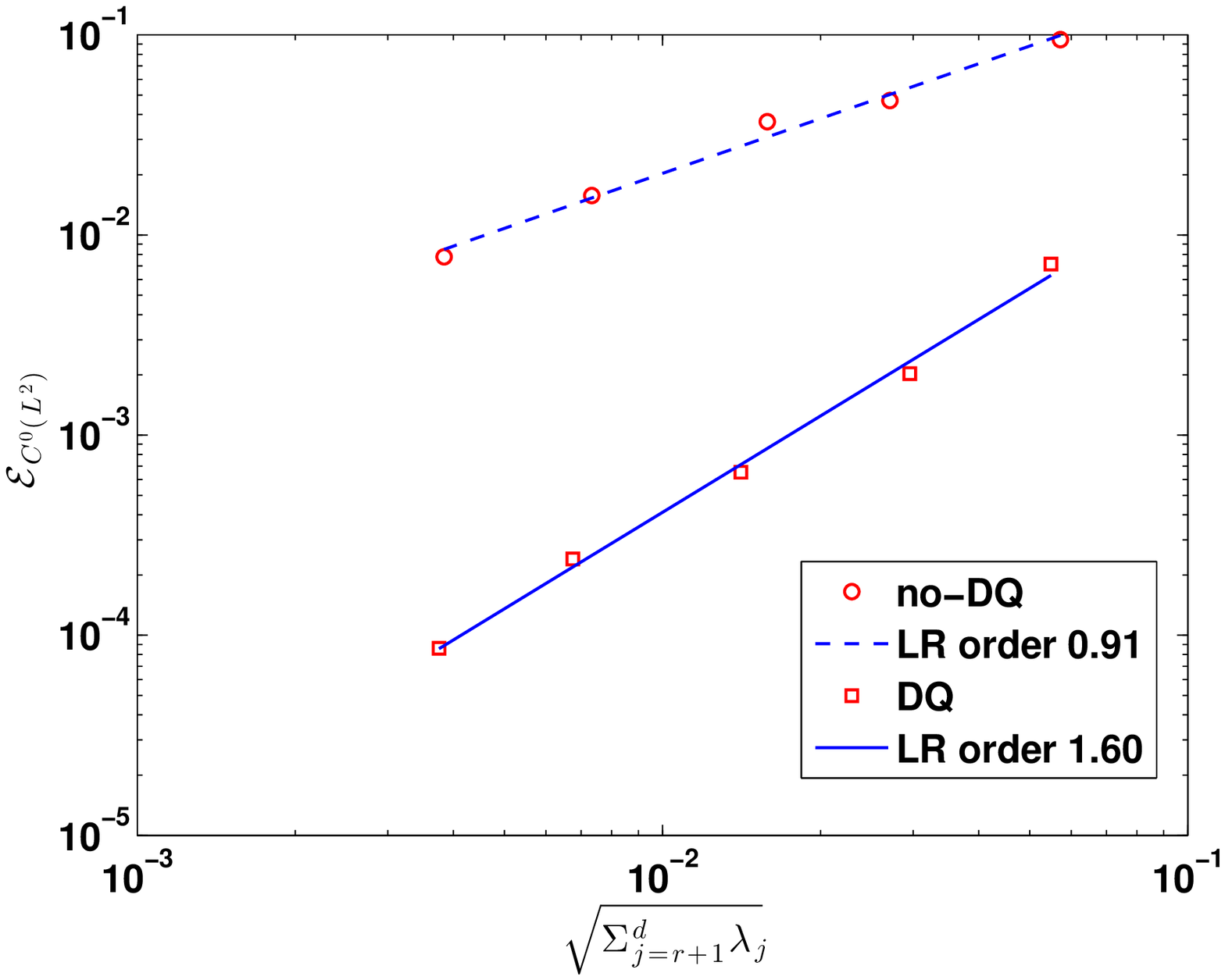}
\end{minipage}
\hspace{.3cm}
\begin{minipage}[h]{.45\linewidth}
\includegraphics[width=1.1\textwidth]{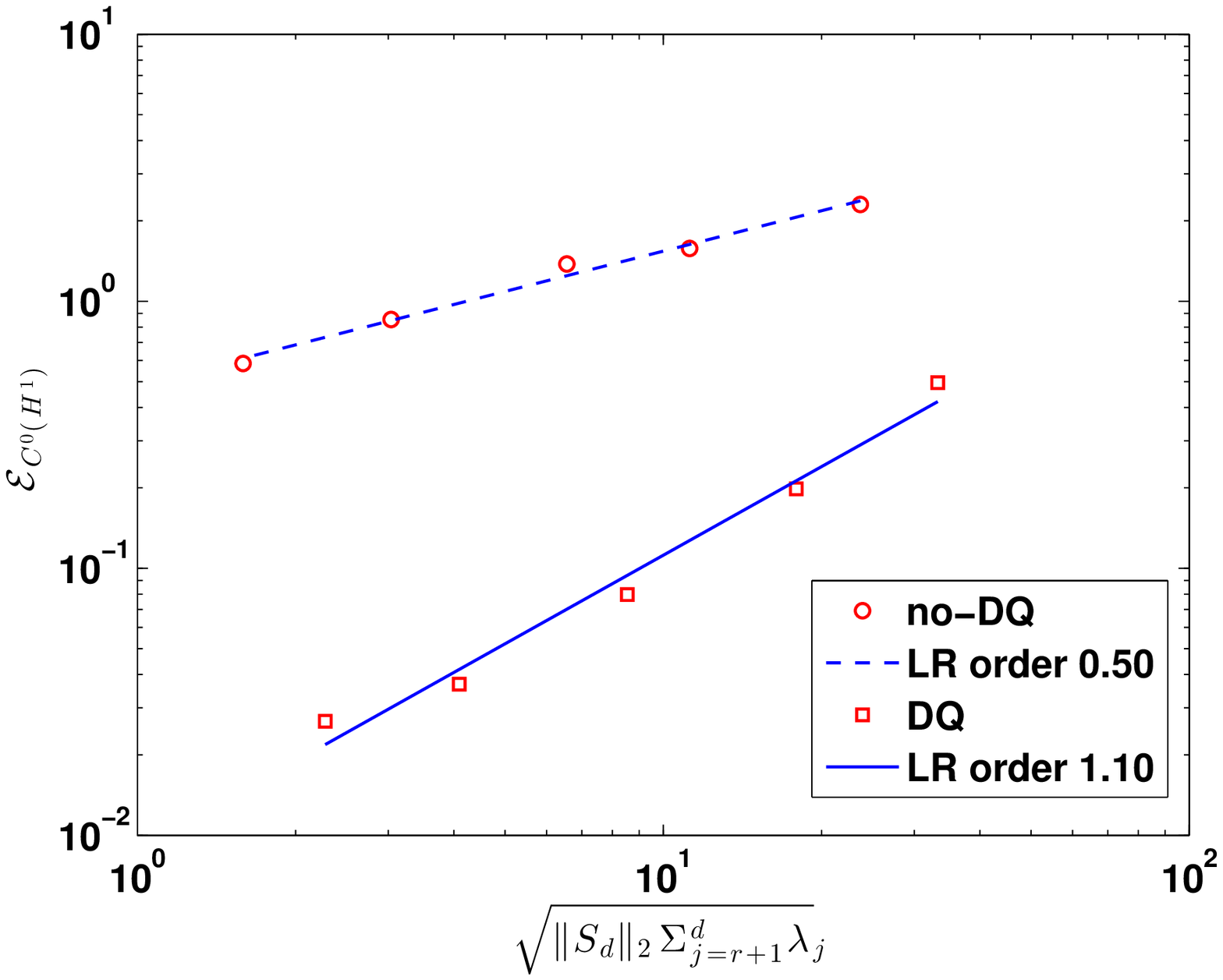}
\end{minipage}
\hspace{.3cm}
\begin{center}
\begin{minipage}[h]{.45\linewidth}
\includegraphics[width=1.1\textwidth]{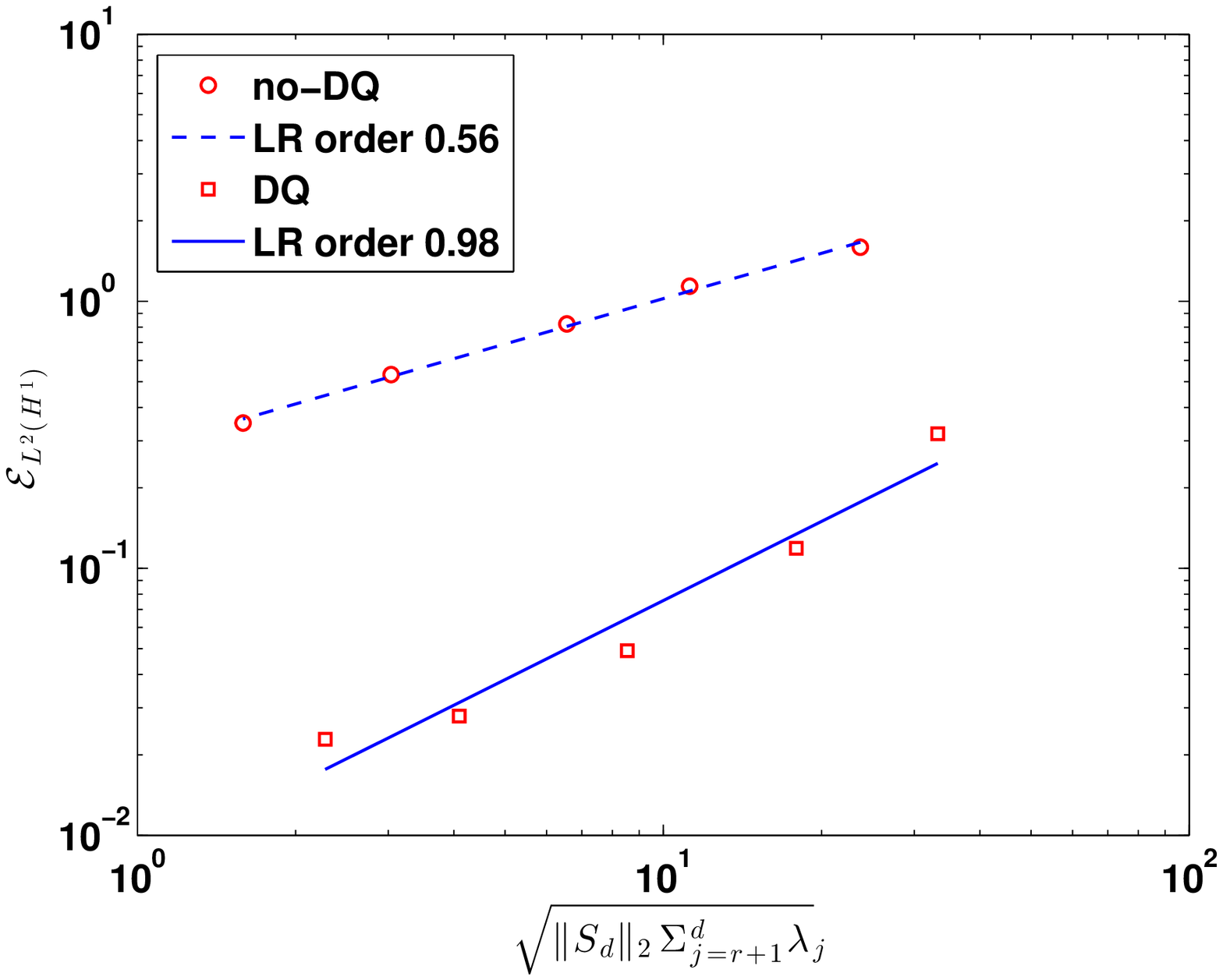}
\end{minipage}
\end{center}
\caption{
	Heat equation.
	Plots of errors in $C^0(L^2)$-norm (top, left), $C^0(H^1)$-norm (top, right), and $L^2(H^1)$-norm (bottom).
}
\label{fig:test1_eC0L2_eC0H1_eL2H1}
\end{figure}

By varying $r$, the number of basis functions used in the POD-G-ROM, we check the rates of convergence with respect to $r$ for the $no\_DQ$ and the $DQ$ cases.
The errors are listed in Table \ref{tab:test1_DQ0} (in the $no\_DQ$ case) and in Table \ref{tab:test1_DQ1} (in the $DQ$ case).
To visualize the rates of convergence with respect to $r$, these errors with their linear regression plots are drawn in Figure \ref{fig:test1_eC0L2_eC0H1_eL2H1}.
The convergence rate of the error in the $C^0(L^2)$-norm, $\cE_{C^0(L^2)}$, is superoptimal in the $DQ$ case and suboptimal in the $no\_DQ$ case.
This supports the theoretical rates of convergence in Table~\ref{tab:numerical_1}, although the suboptimality in the $no\_DQ$ case is mild.
The convergence rate of the error in the $C^0(H^1)$-norm, $\cE_{C^0(H^1)}$, is slightly superoptimal  in the $DQ$ case and strongly suboptimal in the $no\_DQ$ case.
This again supports the theoretical rates of convergence in Table~\ref{tab:numerical_1}.
The convergence rate of the error in the $L^2(H^1)$-norm, $\cE_{L^2(H^1)}$, is optimal in the $DQ$ case and strongly suboptimal in the $no\_DQ$ case.
This supports the theoretical rates of convergence in Table~\ref{tab:numerical_1} for the $DQ$ case, but not for the $no\_DQ$ case.

Overall, the numerical results support the theoretical rates of convergence proved in Section~\ref{sec:error} and summarized in Table~\ref{tab:numerical_1}.
We also emphasize that {\it the convergence rates in the $DQ$ case in all three norms are much higher than (and almost twice as high as) the corresponding rates of convergence in the $no\_DQ$ case}.
\begin{table}[h]
\centering
\tabcaption{Heat equation. Errors in the $no\_DQ$ case.}
\label{tab:test1_DQ0}
\begin{tabular}{|c|c|c|c|c|c}
\hline
 r  & $\Lambda_r$ & $\cE_{C^0(L^2)}$ & $\cE_{C^0(H^1)}$  &  $\cE_{L^2(H^1)}$ \\
 \hline
 3 & 5.72e-02  & 9.46e-02   &   2.30e+00  & 1.59e+00   \\      
 5 & 2.71e-02  & 4.70e-02   &   1.58e+00  & 1.14e+00   \\     
 7 & 1.58e-02  & 3.69e-02   &   1.38e+00  & 8.22e-01     \\    
 10 & 7.34e-03  & 1.57e-02   &   8.54e-01 & 5.31e-01   \\     
 13 & 3.84e-03  & 7.78e-03   &   5.84e-01 & 3.50e-01   \\
 \hline
 \end{tabular}
 \end{table}
\begin{table}[h]
\centering
\tabcaption{Heat equation. Errors in the $DQ$ case.}
\label{tab:test1_DQ1}
\begin{tabular}{|c|c|c|c|c|c}
\hline
 r  & $\Lambda_r$ & $\cE_{C^0(L^2)}$ & $\cE_{C^0(H^1)}$  &  $\cE_{L^2(H^1)}$ \\
 \hline
 19 & 5.49e-02  & 7.15e-03   &   4.96e-01 & 3.19e-01    \\       
 23 & 2.95e-02  & 2.03e-03   &   1.98e-01  & 1.19e-01   \\       
 28 & 1.41e-02  & 6.52e-04   &   7.97e-02  & 4.91e-02   \\      
 33 & 6.75e-03  & 2.41e-04   &   3.68e-02  & 2.80e-02   \\       
 37 & 3.76e-03  & 8.60e-05   &   2.67e-02  & 2.29e-02   \\
 \hline
 \end{tabular}
 \end{table}

\subsection{Burgers Equation}
\label{sec:burgers}

In this section, we consider the one-dimensional Burgers equation.
As mentioned at the beginning of Section~\ref{sec:numerical}, the error estimates proved in Section~\ref{sec:error} are valid for the (linear) heat equation, but not necessarily valid for the nonlinear Burgers equation.
Nevertheless, to gain some insight into the range of validity of the theoretical development in Section~\ref{sec:error}, we investigate the rates of convergence with respect to $r$ in the $no\_DQ$ and the $DQ$ cases for the nonlinear Burgers equation:
 \begin{eqnarray}
\label {burgers}
\left\{ 
\begin{array}{ll}
u_t
- \nu \, u_{xx}
+ u \, u_x
= f
& \qquad \text{ in } \Omega \times (0, T] \, , \\
u(x,0) = u_0(x)
& \qquad \text{ in } \Omega \, , \\
u(x, t) = g(x, t) 
& \qquad \text{ on } \partial \Omega \times (0, T] \, .
\end{array} 
\right.
\end{eqnarray}
The initial condition is
\begin{eqnarray}
u_0(x) 
= \left\{
\begin{array}{cc}
  1 & \quad \text{ if } x \in \left( 0, \frac{1}{2} \right] \\
  \\
  0 & \qquad \text{ if } x \in \left( \frac{1}{2} , 1 \right) \, ,
\end{array}
\right.
\label{ic_1}
\end{eqnarray}
which is similar to that used in \cite{KV01}.
The diffusion parameter is $\nu=10^{-2}$, the forcing term is $f = 0$, $\Om = [0, 1]$, and $T= 1$. 
The boundary conditions are homogeneous Dirichlet, that is, 
$u(0,t) = u(1,t) = 0$ for all $t\in [0, 1]$.  

\begin{figure}[h]
\centering
\begin{minipage}[h]{.6\linewidth}
\includegraphics[width=1\textwidth]{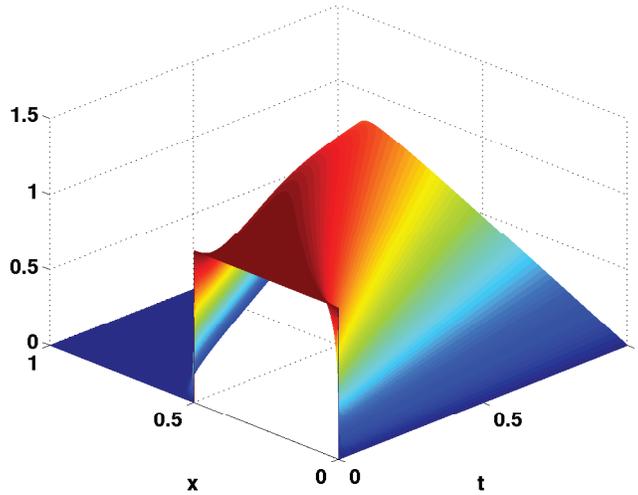}
\end{minipage}
\caption{
	Burgers equation.
	Fine resolution finite element solution used to generate the snapshots.
}
\label{fig:test2dns}
\end{figure}

To generate snapshots, we use piecewise linear finite elements with mesh size $h=1/1024$ and the backward Euler method with a time step $\Delta t = 10^{-4}$, and save data at each time instance. 
The time evolution of the finite element solution is shown in Figure \ref{fig:test2dns}. 
All snapshots are used for the POD basis generation and the same numerical solver is used in the POD-G-ROM. 

\begin{figure}[h]
\begin{minipage}[h]{.45\linewidth}
\includegraphics[width=1.1\textwidth]{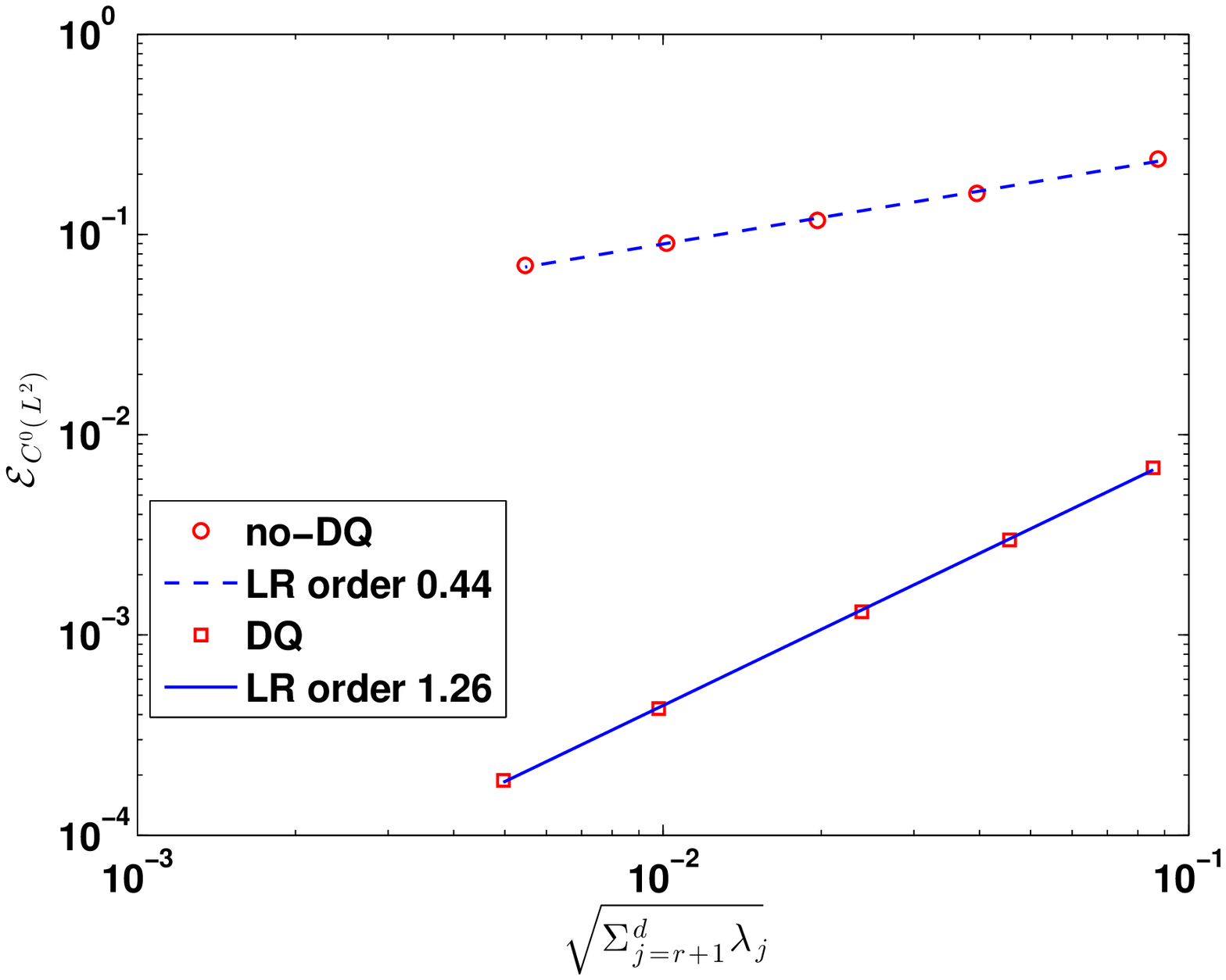}
\end{minipage}
\hspace{.3cm}
\begin{minipage}[h]{.45\linewidth}
\includegraphics[width=1.1\textwidth]{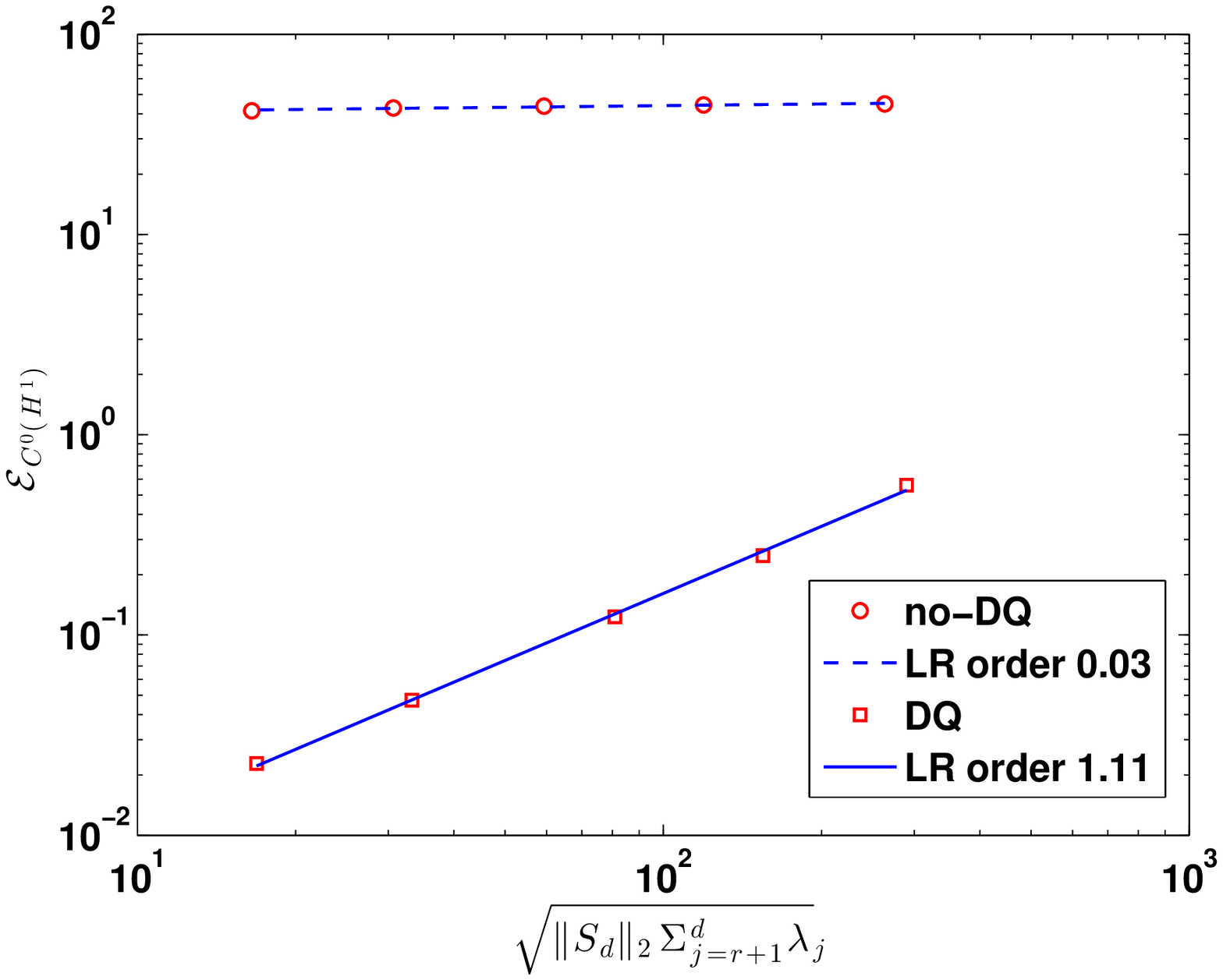}
\end{minipage}
\hspace{.3cm}
\begin{center}
\begin{minipage}[h]{.45\linewidth}
\includegraphics[width=1.1\textwidth]{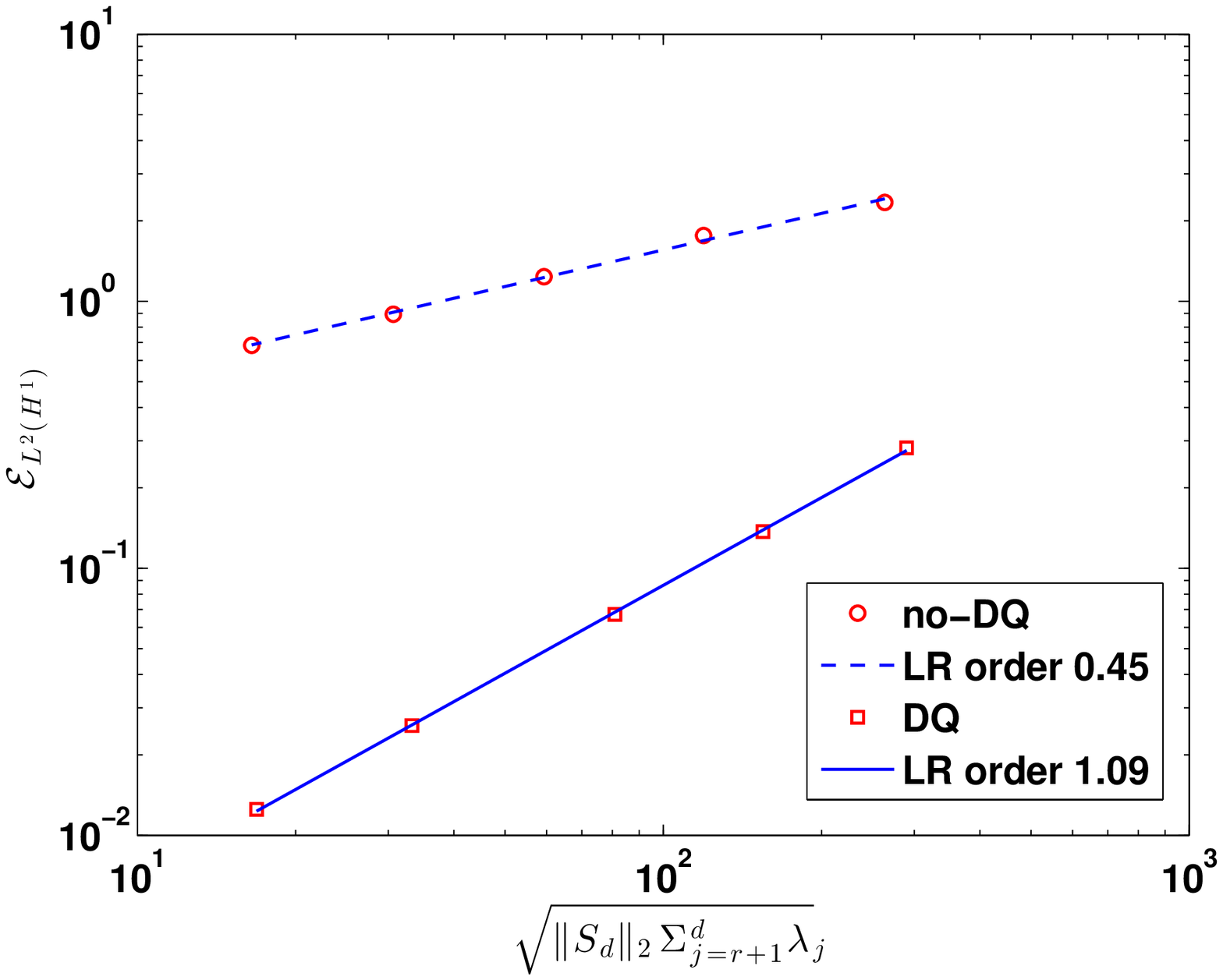}
\end{minipage}
\end{center}
\caption{
	Burgers equation.
	Plots of errors in $C^0(L^2)$-norm (top, left), $C^0(H^1)$-norm (top, right), and $L^2(H^1)$-norm (bottom).
}
\label{fig:test2_eC0L2_eC0H1_eL2H1}
\end{figure}

By varying $r$, we check the rates of convergence with respect to $r$ for the $no\_DQ$ and the $DQ$ cases.
Since the exact solution of the Burgers equation~\eqref{burgers} with the initial condition~\eqref{ic_1} is not known, we consider the errors between the POD-G-ROM results and the snapshots.  
The errors are listed in Table \ref{tab:test2_DQ0} (in the $no\_DQ$ case) and in Table \ref{tab:test2_DQ1} (in the $DQ$ case).
These errors with their linear regression plots are drawn in Figure~\ref{fig:test2_eC0L2_eC0H1_eL2H1}.
The convergence rate of the error in the $C^0(L^2)$-norm, $\cE_{C^0(L^2)}$, is superoptimal in the $DQ$ case and strongly suboptimal in the $no\_DQ$ case.
This clearly supports the theoretical rates of convergence in Table~\ref{tab:numerical_1}.
The convergence rate of the error in the $C^0(H^1)$-norm, $\cE_{C^0(H^1)}$, is superoptimal  in the $DQ$ case and extremely suboptimal in the $no\_DQ$ case.
This strongly supports the theoretical rates of convergence in Table~\ref{tab:numerical_1}.
Finally, the convergence rate of the error in the $L^2(H^1)$-norm, $\cE_{L^2(H^1)}$, is optimal in the $DQ$ case and strongly suboptimal in the $no\_DQ$ case.
This supports the theoretical rates of convergence in Table~\ref{tab:numerical_1} for the $DQ$ case, but not for the $no\_DQ$ case.

Overall, the numerical results clearly support the theoretical rates of convergence proved in Section~\ref{sec:error} and summarized in Table~\ref{tab:numerical_1}.
We also emphasize that {\it the convergence rates in the $DQ$ case in all three norms are much higher than (and at least twice as high as) the corresponding rates of convergence in the $no\_DQ$ case}.
%
\begin{table}[h]
\centering
\tabcaption{Burgers equation. Errors in the $no\_DQ$ approach.}
\label{tab:test2_DQ0}
\begin{tabular}{|c|c|c|c|c|c}
\hline
 r  & $\Lambda_r$ & $\cE_{C^0(L^2)}$ & $\cE_{C^0(H^1)}$  &  $\cE_{L^2(H^1)}$ \\
 \hline
 3 & 8.74e-02  &   2.38e-01  &   4.48e+01  & 2.34e+00   \\       
 5 & 3.95e-02  &   1.60e-01  &   4.44e+01 & 1.76e+00   \\        
 7 & 1.97e-02  &   1.17e-01  &   4.37e+01 & 1.24e+00   \\        
 9 & 1.02e-02  &   9.05e-02  &   4.28e+01  & 8.94e-01    \\       
 11 & 5.47e-03  &   7.01e-02   &   4.14e+01 & 6.84e-01 \\
 \hline
 \end{tabular}
 \end{table}
\begin{table}[h]
\centering
\tabcaption{Burgers equation. Errors in the $DQ$ approach.}
\label{tab:test2_DQ1}
\begin{tabular}{|c|c|c|c|c|c}
\hline
 r  & $\Lambda_r$ & $\cE_{C^0(L^2)}$ & $\cE_{C^0(H^1)}$  &  $\cE_{L^2(H^1)}$ \\
 \hline
 18 & 8.55e-02  &   6.83e-03      &   5.60e-01     & 2.82e-01      \\
 21 & 4.56e-02  &   2.99e-03      &   2.49e-01     & 1.37e-01   \\
 24 & 2.39e-02  &   1.31e-03      &   1.23e-01     & 6.72e-02     \\   
 28 & 9.81e-03  &   4.29e-04      &   4.73e-02     & 2.57e-02    \\    
 31 & 4.97e-03  &   1.88e-04      &   2.28e-02     & 1.25e-02   \\
 \hline
 \end{tabular}
 \end{table}
 
%

\section{Conclusions}
\label{sec:conclusions}

The effect of using or not the snapshot DQs in the generation of the POD basis (the $DQ$ and the $no\_DQ$ cases, respectively) was investigated theoretically and numerically. 
The criterion used in this theoretical and numerical investigation was the rate of convergence with respect to $r$ of the POD-G-ROM solution to the exact solution, where $r$ is the number of POD basis functions used in the POD-G-ROM.

The error estimates in Section~\ref{sec:error} yielded the following conclusions:
In the $DQ$ case, the convergence rates were optimal in all three norms considered (the $C^0(L^2)$-norm, the $C^0(H^1)$-norm and the $L^2(H^1)$-norm.
In the $no\_DQ$ case, the convergence rates were suboptimal in the $C^0(L^2)$-norm and in the $C^0(H^1)$-norm, and optimal in the $L^2(H^1)$-norm.

The numerical results in Section~\ref{sec:numerical} for the (linear) heat equation and the (nonlinear) Burgers equation confirmed the conclusions suggested by the theoretical error estimates in Section~\ref{sec:error}:
In the $DQ$ case, the convergence rates were superoptimal in the $C^0(L^2)$-norm, the $C^0(H^1)$-norm, and the $L^2(H^1)$-norm.
In the $no\_DQ$ case, the convergence rates were suboptimal in the $C^0(L^2)$-norm, the $C^0(H^1)$-norm, and the $L^2(H^1)$-norm.
The only departure from the theoretical conclusions was that, in the $no\_DQ$ case, the convergence rate in the $L^2(H^1)$-norm was suboptimal.  
We emphasize that, for both the heat equation and the Burgers equation,  {\it the convergence rates in the $DQ$ case in all three norms were much higher than (and usually at least twice as high as) the corresponding rates of convergence in the $no\_DQ$ case}.

The theoretical error estimates in Section~\ref{sec:error} and the numerical results in Section~\ref{sec:numerical} strongly suggest the following conjecture:
{\it ``The snapshot DQs should be used in the generation of the POD basis in order to achieve optimal rates of convergence with respect to $r$, the number of POD basis functions utilized in the POD-G-ROM.
}
We also conjecture that using the snapshot DQs in the generation of the POD basis could alleviate some of the degrading of convergence with respect to $r$ seen in, e.g., \cite{rowley2004model,carlberg2011efficient,baiges2012explicit,caiazzo2013numerical}.
We intend to investigate this conjecture in a future study.


\bibliographystyle{plain}
\bibliography{./comprehensive_bibliography}

\end{document}